\documentclass[twoside,tbtags,leqno]{amsart}
\usepackage{amsmath,color,amsthm,amssymb,cite,bbm}

\newcommand{\p}{\mathbb{P}}
\newcommand{\e}{\mathbb{E}}
\newcommand{\va}{\mathbb{V}\mathrm{ar}}

\numberwithin{equation}{section}

\newtheorem{theorem}{Theorem}[section]
\newtheorem{corollary}[theorem]{Corollary}
\newtheorem{lemma}[theorem]{Lemma}
\newtheorem{proposition}[theorem]{Proposition}
\newtheorem{claim}[theorem]{Claim}
\newtheorem{example}[theorem]{\sl Example}
\newtheorem{definition}[theorem]{Definition}

\theoremstyle{definition}
\newtheorem{remark}[theorem]{Remark}

\newcommand{\lzz}{\lim_{z\downarrow 0}}
\newcommand{\lf}{\left\lfloor}

\newcommand{\rf}{\right\rfloor}

\newcommand{\scrs}{{\mathcal{S}}}

\newcommand{\scrb}{\mathcal{B}}

\newcommand{\sep}{{\rm sep}}

\newcommand{\scry}{\mathcal{Y}}

\newcommand{\Ph}{\widehat{P}}

\newcommand{\begp}{\begin{proposition}}
\newcommand{\enp}{\end{proposition}}
\newcommand{\begt}{\begin{theorem}}
\newcommand{\ent}{\end{theorem}}
\newcommand{\begl}{\begin{lemma}}
\newcommand{\enl}{\end{lemma}}
\newcommand{\begc}{\begin{corollary}}
\newcommand{\enc}{\end{corollary}}
\newcommand{\begcl}{\begin{claim}}
\newcommand{\encl}{\end{claim}}
\newcommand{\begr}{\begin{remark}}
\newcommand{\enr}{\end{remark}}
\newcommand{\begal}{\begin{algorithm}}
\newcommand{\enal}{\end{algorithm}}
\newcommand{\begd}{\begin{definition}}
\newcommand{\enf}{\end{definition}}
\newcommand{\begx}{\begin{example}}
\newcommand{\enx}{\end{example}}
\newcommand{\bega}{\begin{array}}
\newcommand{\ena}{\end{array}}

\newcommand{\gD}{\Delta}

\newcommand{\ignore}[1]{}

\newcommand{\sfrac}[2]{{\textstyle\frac{#1}{#2}}}

\def\rompar(#1){\textup(#1\textup)}    

\newcommand\gd{\delta}
\newcommand\gl{\lambda}
\newcommand\gL{\Lambda}

\newcommand\CD{C^{n}}

\newcommand{\refP}[1]{Proposition~\ref{#1}}
\newcommand{\refS}[1]{Section~\ref{#1}}
\newcommand{\refT}[1]{Theorem~\ref{#1}}

\newcommand\ie{i.e.\spacefactor=1000}

\newcommand\essinf{\operatorname{ess\,inf}}

\newcommand\doma{\mathcal{D}_A}
\newcommand\domas{\mathcal{D}_{A^*}}

\newcommand\HD{H^{n}}

\newcommand{\bss}{S^*}

\newcommand\urladdrx[1]{{\urladdr{\def~{{\tiny$\sim$}}#1}}}

\begin{document}
\title[Strong Stationary Duality for Diffusion Processes]
{Strong Stationary Duality for Diffusion Processes}

\author{James Allen Fill}
\address{Department of Applied Mathematics and Statistics,
The Johns Hopkins University,
34th and 
Charles Streets,
Baltimore, MD 21218-2682 USA}
\email{jimfill@jhu.edu}
\urladdrx{http://www.ams.jhu.edu/~fill/}
\thanks{Research supported by the Acheson~J.~Duncan Fund for the Advancement of Research in
Statistics.  
}
\author{Vince Lyzinski}
\address{The Johns Hopkins University Human Language Technology Center of Excellence,
Stieff Building, 810 Wyman Park Drive,
Baltimore, Maryland 21211-2840 USA}
\email{vincelyzinski@gmail.com}
\urladdrx{http://www.ams.jhu.edu/~lyzinski/}
\thanks{Research supported by the Acheson~J.~Duncan Fund for the Advancement of 
Research in Statistics, and by U.S. Department of Education GAANN grant P200A090128.}
\date{}

\begin{abstract}
We develop the theory of strong stationary duality for diffusion processes on compact intervals.  We analytically derive the generator and boundary behavior of the dual process and recover a central tenet of the classical Markov chain theory in the diffusion setting by linking the separation distance in the primal diffusion to the absorption time in the dual diffusion.  We also exhibit our strong stationary dual as the natural limiting process of the strong stationary dual sequence of a well chosen sequence of approximating birth-and-death Markov chains, allowing for simultaneous numerical simulations of our primal and dual diffusion processes.  Lastly, we show how our new definition of diffusion duality allows the spectral theory of cutoff phenomena to extend naturally from birth-and-death Markov chains to the present diffusion context. 
\end{abstract}
\maketitle
\section{Introduction and Background}

Strong stationary duality (SSD)---first developed in the setting of discrete-state Markov chains in \cite{DFSST} and \cite{SSDCT}---has proven to be a powerful tool in the study of mixing times of Markov chains.  In the setting of a suitably monotone Markov chain, strong stationary duality guarantees that a minimal strong stationary time~$T$---whose tail probabilities satisfy $P(T > t) = s(t)$, where $s(t)$ is the separation of the chain at time~$t$---is equal in law to the absorption time in a suitably defined dual chain.  By studying and bounding the absorption time, which is often more tractable than direct consideration of separation, one can tightly bound the separation in the primal chain.  This duality between strong stationary times and absorption times plays a leading role in the development of such diverse techniques as perfect sampling of Markov chains (see \cite{fps}, \cite{MR2001m:60164}), characterizations of separation cut-offs in birth and death chains (see \cite{sepbd}), stochastic constructions of Markov chain hitting times (see \cite{JAFBD}, \cite{PDLmBD}, \cite{LF}), and the analysis of the fastest mixing Markov chain on a graph (see \cite{fk}), to name a few.

However, since initially being referenced in \cite{SSDCT}, extending SSD from Markov chains to the diffusion regime has remained an open problem.  Herein, we present a major step towards this extension.  Utilizing a functional analytic approach, in Section \ref{S:SSD} we systematically develop the theory of SSD for diffusion processes on compact intervals and analytically derive the form of the dual diffusion's generator; in the process, we also explicitly derive the boundary behavior of the dual diffusion.  We further motivate our definition in Section \ref{sec:2} by showing that a suitably defined sequence of Markov chains and their strong stationary duals converge, respectively, to our primal diffusion and its strong stationary dual.  In Section \ref{S:SST}, we recover a central tenet of the classical Markov chain theory in our diffusion setting by linking the separation distance in the primal diffusion to the absorption time in the dual diffusion.  In Section \ref{S:HT}, we exploit this connection to derive the analogue to the birth-and-death cut-off phenomenon theory of \cite{sepbd} in the diffusion setting.

Recently, and independently of our work, wonderful developments in diffusion strong stationary duality have been made in~\cite{mic} and~\cite{mish}.  Our present work was originally presented in the dissertation of the second author \cite{VL} certified in December 2012, and predates the work 
of~\cite{mic} and~\cite{mish}.

\section{Background} 
\label{S:Section 1}  
Let $I=[l,r]\subset\mathbb{R}$ be a compact 
interval (note that the following definitions easily extend to infinite closed $I$).
Denote the interior of $I$ by $I^{\circ}=(l,r)$.  Let $a(\cdot)$, $b(\cdot)\in C(I^{\circ})$ and assume that $b>0$ on $I^{\circ}$.  Then, with appropriate boundary conditions (which we be detailed shortly), the operator
\begin{equation}
\label{eq:gener}
A=\frac{1}{2}b(x)\frac{d^2}{dx^2}+a(x)\frac{d}{dx}.
\end{equation}
acting on $\left\{ f\in C(I)\cap C^2(I^{\circ})\ | \ Af\in C(I)\right\}$
generates a Feller semigroup $(T_t)_{t=0}^{\infty}$ on $C(I)$.  Let $X$ be the time homogeneous one-dimensional diffusion process associated with $A$.

Presently and in the sequel, we will assume $X$ is a regular diffusion process on compact interval $I$ ($ = [0, 1]$, without loss of generality) with initial distribution 
$\Pi_0$ and generator $A$.  Regularity here implies that $\p(T_y<\infty|X_0=x)>0$ for all $x,y\in I^\circ$, where $T_y$ is the hitting time of state $y$ for $X$ (see \cite{kt2} for detail; note that we shall only consider regular diffusions in the sequel). Let $S:I\rightarrow\mathbb{R}$ be a scale function of $X$ defined via
$$
S'(x):=s(x):=\text{exp}\left\{-\int^x\frac{2a(y)}{b(y)}dy\right\},
$$
and let the speed measure of $X$ be denoted by $M$, i.e., $M$ is the nonnegative measure on $I^\circ$ with density $m(x):=[s(x)b(x)]^{-1}.$  Abusing notation, we will also use $M(x)$ to denote the speed function $M(x)=\int_c^x m(y)\,dy$ where $c\in I^\circ$ is fixed but arbitrary.  Feller classified the boundary behavior of $X$ at $l$ (analogous results holding at $r$ if $r\in I$) by looking at the behavior of 
$$N(l):=\int_{(l,x]} \left[S(x)-S(\eta)\right]\,M(d\eta),\ \ \ \Sigma(l):=\int_{(l,x]} \left[M(x)-M(\eta)\right]\, dS(\eta)$$ for a fixed $x\in I^{\circ}$ and  by calculating boundary conditions satisfied by elements of the domain of $A$, which we shall call $\doma$ (see \cite[Section 8.1]{etku} for more details).  Entrance boundaries are characterized by 
$ N(l)<\infty,\ \Sigma(l)=\infty.$  Note that \cite[Section 15.6]{kt2} implies that to show $l$ is entrance, it suffices to show that  $N(l)<\infty$ and $S(l,x]=\lim_{y\downarrow l} [S(x)-S(y)]=\infty$.
Exit boundaries are characterized by 
$N(l)=\infty,\ \Sigma(l)<\infty.$
Natural boundaries are characterized by
$N(l)=\infty$ and $\Sigma(l)=\infty$.  
Finally, regular boundaries are characterized by 
$\ N(l)<\infty,\ \Sigma(l)<\infty.$  The behavior of the diffusion at a regular boundary will be characterized by boundary conditions satisfied by elements $f\in\doma$.  In particular, we say that $l$ is \emph{instantaneously reflecting} if $f\in\doma$ implies that 
$$\frac{df}{dS}^+(l)=\lim_{x\downarrow l} \frac{df}{dS}(x)=\lim_{x\downarrow l} \frac{f'(x)}{s(x)}=0.$$  We say $l$ is \emph{absorbing} if $f\in\doma$ implies that $(Af)(l)=0$. 

Assume that~$0$ 
and~$1$ are instantaneously reflecting boundaries for $X$ (and note that analogous results and definitions will hold in the entrance boundary case). The boundary behavior of $X$ guarantees that $M$ is a finite measure on $I^{\circ}$, and normalizing $M(dx)$ to a probability measure gives the unique invariant  distribution of $X$, which we will denote  by $\Pi(dx)$.  For arbitrary $c\in I^{\circ}$, let  us adopt the shorthand $\Pi(x):=\int_{y=c}^x\pi(y)\,dy$, where $\pi$ is the density for $\Pi$ with respect to Lebesgue measure, and note that regularity of $X$ guarantees $\pi>0$ on $I^{\circ}$.  The reflecting behavior at $0$ guarantees $\lim_{c\downarrow 0} \int_{y=c}^x\pi(y)\,dy$ exists and is finite for all $x\in I^{\circ}$, and so to ease notation we may let $\Pi(x)=\int_{y=0}^x\pi(y)\,dy$ defined as an improper integral.  Lastly, let $(P_t)_{t=0}^{\infty}$ be the Markov transition function associated with $X$ and denote the corresponding transition densities with respect to Lebesgue measure by $(p_t)_{t=0}^{\infty}$.

Based on the boundary behavior of $X$, we can completely specify the domain of $A$ as 
\begin{equation}
\label{eq:gen}
\mathcal{D}_A= \left\{ f\in C(I)\cap C^2(I^{\circ})\ | \ Af\in C(I),\ \frac{df}{dS}^+(0)=\frac{df}{dS}^-(1)=0 \right\}
\end{equation}
 (see \cite[Section~8.1]{etku}, especially~(1.11) there, with $q_0 = 0 = q_1$ because both boundaries are instantaneously reflecting), 
where as above
$$ \frac{df}{dS}^+(0)=\lim_{x\downarrow 0} \frac{df}{dS}(x)=\lim_{x\downarrow 0} \frac{f'(x)}{s(x)},$$ and
$$ \frac{df}{dS}^-(1)=\lim_{x\uparrow 1} \frac{df}{dS}(x)=\lim_{x\uparrow 1} \frac{f'(x)}{s(x)}.$$

Let 
$\mathcal{F}[0,1]$ be the space of bounded real valued measurable functions on $[0,1]$ equipped with its usual Borel 
$\sigma$-field $\mathcal{B}$.  Let $\mathcal{M}[0,1]$ be the space of signed measures on $([0,1], \mathcal{B}).$  As  
in~\cite[Section 7.1]{lamp}, we note the natural bilinear functional on $\mathcal{F}[0,1] \times \mathcal{M}[0,1]$ defined by 
$(\mu,f)=\int_0^1f(x)\,\mu(dx)$.  We denote the adjoint of the operator $T_t$ (with respect to this functional) by $U_t$, where $(T_t)_{t=0}^{\infty}$ is the one parameter Markov semigroup associated 
with $(P_t)_{t=0}^{\infty}$.  

Note that $a(\cdot),\ b(\cdot),$ $m(\cdot),$ and $\pi(\cdot)$ are defined only on $I^{\circ}$.  For notational convenience, any expressions involving these functions and $\partial I$ are to be interpreted as the corresponding limiting expression (when such a limit exists!).  For example, for $0<x<1$ we shall write the improper integral $\int_0^x f(y)\pi(y)\,dy$ rather than the equivalent $\lim_{z\downarrow 0} \int_z^x\!f(y)\pi(y)\,dy$.
\section{Strong stationary duality for diffusions}
\subsection{Definition of the strong stationary dual}
\label{S:SSD}
Let $X^*$ be a second (Feller) diffusion process on $I$ with initial distribution $\Pi_0^*$ and generator $A^*$.
As in the continuous-time discrete-state Markov chain setting (see \cite{SSDCT}), we define the notion of algebraic duality between~$X$ and $X^*$:
\begd
\label{D:1.1}
Consider the integral operator~$\Lambda$ acting on $F[0, 1]$ defined by
$$
(\Lambda f)(x) := 
\begin{cases}
\int_0^x\!\pi^{(x)}(y)f(y)\,dy & \mbox{\ if $x>0$},\\
f(0) & \mbox{\ if $x = 0$},
\end{cases}
$$
where we define the kernel
$$
\pi^{(x)}(y) := \frac{\pi(y)}{\Pi(x)}\ \ \  0 < y \leq x<1,\ \ \text{and }\pi^{(1)}\equiv\pi. 
$$
We say that $X^*$ is a strong stationary dual of~$X$ if 
\begin{equation}\label{eq:dom} 
\mbox{$\gL$ maps $\mathcal{D}_A$ into $\mathcal{D}_{A^*}$}  
\end{equation}
and
\begin{equation} \label{eq:link}
\Lambda A=A^*\Lambda\  \text{as operators defined on $\mathcal{D}_A$}  
\end{equation}
and
\begin{equation}\label{eq:init} 
(\Pi_0,f)=(\Pi_0^* ,\gL f)\ \text{for all f}\in \mathcal{F}[0,1].
\end{equation}  
\enf
\begr
\label{rem:contgl}
If $f\in C(I)$, then $\gL f\in C(I)$ as well.  To show this, first note that $\pi\in C(I^{\circ})$, $\Pi\in C(I)$, and for $x>0$ we have $\Pi(x)>0$.  Clearly, then, $$\gL f(x)=\frac{\int_0^x\! \pi(y) f(y)\,dy}{\Pi(x)}$$ is continuous at all $x>0$.  Continuity at zero is immediate as for any $\epsilon>0$, we can choose $x$ such that $|f(y)-f(0)|<\epsilon$ for all $y\leq x$, and so
$$|\gL f(0)-\gL f(x)|=\left|\int_0^x\! (f(0)-f(y))\pi^{(x)}(y)\,dy\right|\leq \epsilon.$$
\enr
\begr
\label{rem:init}
For $x<1$, let  $\Pi^{(x)}$ be the distribution~$\Pi$ conditioned to $(0, x]$, so that $\Pi^{(x)}$ 
has density $\pi^{(x)}$ when $x > 0$, and let $\Pi^{(0)} := \delta_0$ and $\Pi^{(1)}:=\Pi$.  If $\Pi_0=\Pi^{(x)}$ for some $x\in[0,1]$, then \eqref{eq:init} is uniquely satisfied 
by $\Pi_0^*=\delta_x$.  For $x\in(0,1)$, this is easily seen via
\begin{equation}\label{eq:rem1}
\begin{split}
 \int_I \!f(y)\,\pi^{(x)}(y)\,dy&=(\Pi_0,f)\\
&=(\Pi_0^*,\gL f)=\int_I \int_{y\in(0,z]}\!\pi^{(z)}(y)f(y)\,dy\,\Pi_0^*(dz)\\
&=\int_I \int_{z\in[y,1]}\!\pi^{(z)}(y)\,\Pi_0^*(dz)\,f(y)\,dy. 
\end{split}
\end{equation}
Letting $f(y)=\mathbbm{1}(y>x)$ we see $\Pi^*_0$ must be concentrated on $(0, x]$.
It also follows that for almost every~$y$ satisfying $0<y\leq x$ we have
$$
\pi^{(x)}(y)=\int_{z\in[y,1]}\!\pi^{(z)}(y)\,\Pi_0^*(dz)=\int_{z\in[y,x]}\!\pi^{(z)}(y)\,\Pi_0^*(dz),
$$ 
or, equivalently,
$$
\frac{1}{\Pi(x)} = \int_{z \in [y, x]} \!\frac{\Pi^*_0(dz)}{\Pi(z)}.
$$ 
Letting $y \uparrow x$ through such values,
it follows that $\Pi^*_0=\delta_x$ is the only possible initial distribution for $X^*$.  To show that $\Pi^*_0=\delta_x$ satisfies (\ref{eq:init}), note 
$$
(\delta_x,\gL f)=\gL f(x)=\int_0^x\! \pi^{(x)}(y)f(y)\,dy=(\Pi^{(x)},f),
$$ as desired.
For $x = 0$, the argument goes as follows.  For uniqueness,
if $\Pi_0=\delta_0$, then letting $f(y) = \mathbbm{1}(y \in (0,1])$, the left side 
of~\eqref{eq:init} equals $f(0) = 0$, and the right side is strictly positive unless $\Pi_0^* = \delta_0$.   To see that $\Pi^*_0 = \delta_0$ satisfies~\eqref{eq:init} when $\Pi_0 = \Pi^{(0)} = \delta_0$, we compute 
$(\delta_0,\gL f)=\gL f(0)=f(0)=(\delta_0,f)$.  For $x=1$, the argument proceeds as follows.  Let $f(y)=\mathbbm{1}(y\in [1/2,1])$.  Then 
$$\Lambda f(z)=\begin{cases}
1-\frac{\Pi(1/2)}{\Pi(z)}&\text{ if }z>1/2,\\
0&\text{ otherwise }.
\end{cases}
$$  As $\pi>0$ on $I^{\circ}$, it follows that $\Pi$ is strictly increasing on $I^\circ$ and so 
$\Lambda f(z)<(\Pi,f)$ for all $z<1$.  This implies that equation \eqref{eq:init} is uniquely satisfied by $\Pi_0^*=\delta_1$ as desired.
\enr

\subsection{The dual generator}
\label{S:DG}
From the definition of strong stationary duality, we derive the form of the dual generator:
\begt
\label{theorem:dual}
Let $X$ be defined as above and assume that $X$ has instantaneously reflecting boundaries at $0$ and $1$.  Assume further that 
$b\in C^1(I^{\circ})$.
If $X^*$ is 
a strong stationary dual of~$X$, then the generator $A^*$  of $X^*$ has the form 
$$
(A^*f)(x) = \left(\frac{1}{2}b'(x)-a(x)+b(x)\frac{\pi(x)}{\Pi(x)}\right)f'(x)+\frac{1}{2}b(x)f''(x)
$$
for $x\in I^{\circ}$ and $f\in\mathcal{D}_{A^*}$.  Also
$0$ is an entrance boundary for $X^*$ and 1 is a regular absorbing boundary of $X^*$.
\ent
\begin{remark}
\label{rem:gen}
The dual generator's domain then can be explicitly defined as
$$\mathcal{D}_A^*= \left\{ f\in C(I)\cap C^2(I^{\circ})\ | \ A^*f\in C(I),\ A^*f(1-)=0 \right\}
$$
\end{remark}
\begin{proof}
Let $f\in\doma$.  Then $Af\in C(I)$ and for $x>0$ we have
\begin{align*}
(\gL Af)(x)&=\int_0^x\!\left( a(y)f'(y)\pi^{(x)}(y)+\sfrac{1}{2}b(y)f''(y)\pi^{(x)}(y) \right)\,dy\\
&=\frac{1}{\Pi(x)} \int_0^x\!\left( a(y)f'(y)\pi(y)+\sfrac{1}{2}b(y)f''(y)\pi(y)\right)\,dy.
\end{align*}
We know that there exists a nonzero constant $C$ such that $C\cdot\pi(x)=m(x)$, so that $\pi(x)=\frac{1}{Cb(x)s(x)}$.  Also, $\frac{d}{dx}\frac{1}{s(x)}=\frac{1}{s(x)}\frac{2a(x)}{b(x)}.$ 
We can then rewrite 
\begin{align*}(\Lambda Af)(x)=&\frac{1}{\Pi(x)}\frac{1}{2C}\int^x_0\!\left( \frac{2a(y)}{b(y)}\frac{f'(y)}{s(y)}+\frac{f''(y)}{s(y)} \right)\,dy\\
=&\frac{1}{\Pi(x)}\frac{1}{2C}\int^x_0 \frac{d}{dy} \frac{f'(y)}{s(y)}\,dy\\
=&\frac{1}{\Pi(x)}\frac{1}{2C}\left[ \frac{df}{dS}(x)-\frac{df}{dS}^+(0) \right]\\
=& \frac{1}{2}\frac{b(x)\pi(x)}{\Pi(x)}f'(x),
\end{align*}
as $0$ being a reflecting boundary of~$X$ and $f \in \mathcal{D}_A$ implies $\frac{df}{dS}^+(0)=0$.
 
Let $g\in\domas$.  
For $x\in(0,1)$, from equation (\ref{eq:gener}) for $A^*$ we can write 
$$
(A^*g)(x)=a^*(x)g'(x)+\sfrac{1}{2}b^*(x)g''(x).
$$ 
for some $a^*$, $b^*\in C(I^{\circ}).$  If $f\in\doma$ then by~\eqref{eq:dom} we have $\gL f\in\domas$, and so for $x\in(0,1)$ we know
$(A^*\gL f)(x)=a^*(x)(\gL f)'(x)+\frac{1}{2}b^*(x)(\gL f)''(x)$.  Note that $Af\in C(I)$ by assumption and so $\gL Af=A^*\gL f\in C(I)$ from Remark \ref{rem:contgl}.  Now
\begin{align*}
(\gL f)'(x)&=\frac{\Pi(x)\pi(x)f(x)-\pi(x)\int_0^x\pi(y)f(y)\,dy}{\Pi(x)^2}\\
&=\frac{\pi(x)}{\Pi(x)}[f(x)- (\gL f)(x)]
\end{align*}
and so
\begin{align*}
(\gL f)''(x)&=\frac{\Pi(x)\pi'(x)-\pi(x)^2}{\Pi(x)^2}[f(x)- (\gL f)(x)]+\frac{\pi(x)}{\Pi(x)}\left\{f'(x)-\frac{\pi(x)}{\Pi(x)}[f(x)
- (\gL f)(x)]\right\}\\
&=\left[\frac{\pi'(x)}{\Pi(x)}-\frac{2\pi(x)^2}{\Pi(x)^2}\right][f(x)- (\gL f)x)]+\frac{\pi(x)}{\Pi(x)}f'(x).
\end{align*}
Now 
by~\eqref{eq:link}, $\Lambda A=A^*\Lambda$ as operators on~$\doma$, which implies that for any $x\in(0,1)$ and $f\in\doma$
we have
\begin{align}
\label{eq:linked}
\frac{1}{2}\frac{b(x)\pi(x)}{\Pi(x)}f'(x)=\bigg(&a^*(x)\frac{\pi(x)}{\Pi(x)}+\frac{1}{2}b^*(x)\bigg[\frac{\pi'(x)}{\Pi(x)}-\frac{2\pi(x)^2}{\Pi(x)^2}\bigg]\bigg)[f(x)- (\gL f)(x)]  \nonumber \\
&+\frac{1}{2}b^*(x)\frac{\pi(x)}{\Pi(x)}f'(x). 
\end{align}
For any fixed $x\in I^{\circ}$, we can 
choose $f\in\doma$ so that $f'(x)=0$ and $f(x)\neq (\gL f)(x)$ [e.g., let $f$ be a suitably smooth approximation of 
$\mathbbm{1}(x/3,x/2)$], and for any such~$f$, equation~\eqref{eq:linked} yields 
\begin{equation}\label{eq:bstar}
a^*(x)\frac{\pi(x)}{\Pi(x)}+\frac{1}{2}b^*(x)\bigg[\frac{\pi'(x)}{\Pi(x)}-\frac{2\pi(x)^2}{\Pi(x)^2}\bigg]=0.
\end{equation}
We then find for $f \in \doma$ and $x \in (0, 1)$ that 
$(A^*\gL f)(x)=\frac{1}{2}\frac{b^*(x)\pi(x)}{\Pi(x)}f'(x)$,
and by~\eqref{eq:link} this equals $(\gL A f)(x) = \frac{1}{2}\frac{b(x)\pi(x)}{\Pi(x)}f'(x)$.  For each $x$ in $(0,1)$, we can choose an $f\in\doma$ such that $f'(x)\neq 0$, and using any such~$f$ we find that $b^*(x) = b(x)$. 

Next, we have from $\pi(x)=\frac{1}{Cb(x)s(x)}$ that 
$$
\pi'(x)=\frac{-b'(x)s(x)-b(x)s'(x)}{Cb(x)^2s(x)^2}.
$$
Equation~\eqref{eq:bstar} and $b^* \equiv b$ then yields  
\begin{align*}
\frac{\pi(x)}{\Pi(x)}a^*(x)&=\frac{1}{2}b(x)\left[\frac{b'(x)s(x)+b(x)s'(x)}{C\Pi(x)b(x)^2s(x)^2}+\frac{2\pi(x)^2}{\Pi(x)^2}\right]\\
&=\frac{1}{2C\Pi(x)}\left[\frac{b'(x)}{b(x)s(x)}+\frac{s'(x)}{s(x)^2}\right]+b(x)\frac{\pi(x)^2}{\Pi(x)^2}\\
&=\frac{1}{2C\Pi(x)}\left[Cb'(x)\pi(x)-\frac{2a(x)}{s(x)b(x)}\right]+b(x)\frac{\pi(x)^2}{\Pi(x)^2}\\
&=\frac{1}{2}b'(x)\frac{\pi(x)}{\Pi(x)}-a(x)\frac{\pi(x)}{\Pi(x)}+b(x)\frac{\pi(x)^2}{\Pi(x)^2},
\end{align*}
so that $a^*(x)=\frac{1}{2}b'(x)-a(x)+b(x)\frac{\pi(x)}{\Pi(x)}$ on $I^\circ$, as desired.

To find the boundary behavior of the dual diffusion at~$0$ and at~$1$, we calculate the dual scale function and the dual speed measure.  First, note that
\begin{align}
s^*(x)&=\exp\left[ -\int^x\!\frac{2a^*(y)}{b^*(y)}\,dy \right] \nonumber \\
&=\exp\left[ -\int^x\!\frac{b'(y)}{b(y)}\,dy+\int^x\!\frac{2a(y)}{b(y)}\,dy-\int^x\!\frac{2m(y)}{M(y)}\,dy \right] \nonumber \\
&=\frac{1}{b(x)}\frac{1}{s(x)}\frac{1}{M(x)^2} \nonumber \\
\label{eq:dualscale}
&=\frac{m(x)}{M(x)^2}, 
\end{align} 
and a scale function for $X^*$ is
\begin{equation}
\label{eq:dualScale}
S^*(x)=\frac{-1}{M(x)}.
\end{equation}
  
Next, note 
\begin{equation}
\label{eq:dualspeed}
m^*(x)=\frac{1}{b^*(x)s^*(x)}=\frac{M(x)^2}{m(x)b(x)}=M(x)^2s(x).
\end{equation} 
Now $M(x)$ is continuous on $I$ and 
$M(0)=0$, so
there is a~$y$ such that $M(\zeta) \leq 1$ for all $\zeta \leq y$.  For the dual scale measure $S^*$ we then have 
\begin{align*}
S^*(0,y]&= \int _{(0,y]}\! s^*(\zeta)\,d\zeta\\
&= \lzz \int_{z}^y \frac{m(\zeta)}{M(\zeta)^2}\,d\zeta\\
&\geq \lzz \int_{z}^y \frac{m(\zeta)}{M(\zeta)}\,d\zeta\\
&= \lzz \big[ \log M(y) - \log M(z) \big] = \infty.
\end{align*}  
To show that~$0$ is an entrance boundary for $X^*$, 
it now suffices to show that $N^*(0)<\infty$.  This is shown via
\begin{align*}
N^*(0)
&=\lzz \int_z^x\!S^*[y,x]\,dM^*(y)\\
&= \lzz \int_z^x\![S^*(x)-S^*(y)]\,m^*(y)\,dy\\
&= \lzz \int_z^x\!\left[ \frac{-1}{M(x)}-\frac{-1}{M(y)} \right] M(y)^2\,s(y)\,dy\\
&\leq \frac{-1}{M(x)} \liminf_{z\downarrow 0} \int_z^x\!M(y)^2\,s(y)\,dy + \limsup_{z\downarrow 0}\int_z^x\!M(y)\,s(y)\,dy.
\end{align*}
It now clearly suffices to prove
$\int_0^x\!M(y)\,s(y)\,dy < \infty$, which follows from the following calculation:
$$
\int_0^x\!M(y)\,s(y)\,dy = \int_0^x\!M(y)\,dS(y) = \int_0^x\,S[y,x]\,dM(y) =: N(0)<\infty,
$$
where we used the fact that $0$ is a reflecting boundary for~$X$ and hence $S(0,x]<\infty$ (else $0$ would be entrance) to derive the second equality and the final inequality.

To prove that~$1$ is a regular absorbing boundary for $X^*$, we first show that  $N^*(1)<\infty.$  Indeed, for fixed $x$ in $I^{\circ}$ we have 
[using \eqref{eq:dualScale}--\eqref{eq:dualspeed}] that
\begin{align*}
N^*(1)&=\int_{[x,1)}\![S^*(y)-S^*(x)]m^*(y)\,dy=\int_{[x,1)}\!\left[\frac{1}{M(x)}-\frac{1}{M(y)}\right]s(y)M^2(y)\,dy\\
&=\int_{[x,1)}\!s(y) \frac{M^2(y)}{M(x)}\,dy-\int_{[x,1)}\!s(y) M(y)\,dy\\
&\leq \frac{M(1)}{M(x)}\int_{[x,1)}\!s(y)M(y)\,dy-\int_{[x,1)}\!s(y) M(y)\,dy<\infty
\end{align*}
where the finiteness holds since~$1$ is reflecting for $X$ [hence $\Sigma(1)<\infty$] and $M(\cdot)$ is increasing and bounded on $I^{\circ}$.  
The finiteness of $N^*(1)$ implies that $1$ is either an entrance or regular boundary for $X^*$.

To show that the boundary is regular absorbing, it suffices to prove that $P^*_1(X^*_t=1)=1$ (as $1$ is then not an entrance boundary by definition, and by the infinitesimal characterization of the generator $A$, we immediately have that that $A^*f(1)=0$ for all $f$ in $\domas$ and hence $1$ is absorbing), for which we will use \refP{prop:sgl} below.  
From that proposition, for any $f \in C[0, 1]$ and $x \in [0, 1]$, we 
have $(\gL T_t f)(x) =(T_t^*\gL f)(x)$.  When $x > 0$, we then have
\begin{align*}
\int_0^1\left[ \int_0^x \pi^{(x)}(y) p_t(y,z)\,dy \right] f(z)\,dz
&= \int_{[0, 1]}\!\int_0^y\!\pi^{(y)}(z) f(z)\,dz\,P^*_x\,(X^*_t\in dy) \\
&= \int_0^1 \left[ \int_{[z, 1]}\!\pi^{(y)}(z) P^*_x\,(X^*_t \in dy) \right] f(z)\,dz.
\end{align*}
In particular, letting $x = 1$ we find
$$
\int_0^1\!\,\pi(z)\,f(z)\,dz
= \int_0^1 \left[ \int_{[z, 1]}\!\pi^{(y)}(z)\, P^*_1(X^*_t \in dy) \right] f(z)\,dz.
$$
Since this holds for all $f \in C[0,1]$, and since both $\pi(z)$ and the expression in square brackets on the right are continuous functions of~$z \in (0, 1]$, it follows, for all $z \in (0, 1]$, that 
$$
\pi(z) = \int_{[z, 1]} \frac{\pi(z)}{\Pi(y)}\,P^*_1(X^*_t\in dy),
$$
and hence $\int_{[z,1]} \frac{1}{\Pi(y)}\,P^*_1(X^*_t\in dy)=1$.  It now follows that
$P^*_1(X^*_t=1)=1$ as desired.
\end{proof}

\begp
\label{prop:sgl}
Let $X^*$ be a strong stationary dual of~$X$, and let the one-parameter Markov semigroups of operators for $X^*$ and~$X$ be $(T^*_t)$ and $(T_t)$ respectively.  Then for all~$t$ we
have $\gL T_t = T_t^*\gL$ as operators on $C[0,1]$.
\enp
\begin{proof}
For all $\lambda$ we have $\gL(\lambda I-A)=(\lambda I-A^*)\gL$ and so the resolvent operators satisfy 
$\gL R_{\lambda}=R_{\lambda}^* \gL$.  For $f \in C[0,1]$ and $x \in [0, 1]$, note that 
$$
(R_{\lambda}^*\gL f)(x) = \int_0^{\infty}\!e^{-\lambda t} (T_t^*\gL f)(x)\,dt
$$
and that 
\begin{align*}
(\gL R_{\gl} f)(x)
&=\int^x_0\!\pi^{(x)}(y)\,(R_{\gl}f)(y)\,dy\\
&=\int_0^x\!\int_0^{\infty}\!\pi^{(x)}(y)\,e^{-\gl t}\,(T_t f)(y)\,dt\,dy\\
&=\int_0^{\infty}\!e^{-\lambda t}\,(\gL T_t f)(x)\,dt.
\end{align*}
Now, by the uniqueness of Laplace transforms of real valued bounded functions, we have $(\gL T_t f)(x) = (T_t^*\gL f)(x)$
for almost all $t$.  Extension to all $t$ follows from the continuity in $t\geq 0$ of $(T_t f)(x)$ [resp.,\ $T^*_t g(x)$] for all $f\in\doma$ (resp.,\ $g$ in $\domas$); see \cite[Chapter 7]{lamp}.
\end{proof}
\begr
\label{rem:sgmble}
From \refP{prop:sgl}, we have that $\Lambda T_t = T^*_t \Lambda$ as operators on $C[0, 1]$ which implies that the equality also holds as operators on $\mathcal{F}[0, 1]$.
\enr

The choice of~$0$ and~$1$ as instantaneously reflecting boundaries was done to streamline exposition.  However, we can establish analogues of ~\refT{theorem:dual} for more general boundary behaviors of $X$.  If 0 and 1 are entrance boundaries for $X$, then the domain of $A$ is 
$$\mathcal{D}_A=\{f\in C(I)\cap C^2(I^{\circ})\ | \ Af\in C(I)\}.$$  If 0 (resp.,\ 1)  is made reflecting then 
we impose the extra condition that $\frac{df}{dS}^+(0)=0$ [resp.,\ $\frac{df}{dS}^-(1)=0$] for functions $f\in\mathcal{D}_A$.
In the proof of ~\refT{theorem:dual}, only the following properties of the boundary at 0 were needed:
$$\frac{df}{dS}^+(0)=0\  \mathrm{for}\ f \in C(I),\ \ N(0)<\infty,$$
and these properties also hold if $0$ is an entrance boundary; 
see \cite[Theorem 12.2]{kt2} for proof that $\frac{df}{dS}^+(0)=0\  \mathrm{for}\ f \in C(I)$ in the entrance boundary case.

Absorption of $X^*$ at 1 is proven completely analogously to the reflecting case.  If $1$ is an entrance boundary for~$X$, then $1$ is an exit boundary for $X^*$ since
\begin{align*}
N^*(1)&=\int_{[x,1)}\![S^*(y)-S^*(x)]m^*(y)\,dy=\int_{[x,1)}\!\left[\frac{1}{M(x)}-\frac{1}{M(y)}\right]s(y)M^2(y)\,dy\\
&=\int_{[x,1)}\!s(y) \left[\frac{M^2(y)}{M(x)}- M(y)\right]\,dy\\
&\geq \int_x^1\!s(y)[M(y)- M(x)]\,dy=\Sigma(1)=\infty
\end{align*}
and integration by parts yields
\begin{align*}
\Sigma^*(1)&=\lim_{z\uparrow 1}\int_x^z\![S^*(z)-S^*(y)]m^*(y)\,dy\\
&=\lim_{z\uparrow 1}\int_x^z\!\left[\frac{1}{M(y)}-\frac{1}{M(z)}\right]s(y)M^2(y)\,dy\\
&=\lim_{z\uparrow 1}\int_x^z\!s(y)\left[ M(y)-\frac{M^2(y)}{M(z)}\right]\,dy\\
&\leq \lim_{z\uparrow 1}\int_x^z\!s(y)[M(z)-M(y)]\,dy\\
&= \lim_{z\uparrow 1} \left( [S(z)-S(x)][M(z)-M(z)]+\int_x^z[S(y)-S(x)]m(y)\,dy \right)\\
&=N(1)<\infty.
\end{align*}
We thus arrive at the following generalization of \refT{theorem:dual}.
\begt
\label{thm:entrefl}
Let $X$ be a regular diffusion on $I$, and assume that each of the boundary points of $I$ is either reflecting or entrance.  Assume further that $b\in C^1(I^{\circ})$.
If $X^*$ is 
a strong stationary dual of~$X$, then the generator $A^*$  of $X^*$ has the form 
$$
(A^*f)(x) = \left(\frac{1}{2}b'(x)-a(x)+b(x)\frac{\pi(x)}{\Pi(x)}\right)f'(x)+\frac{1}{2}b(x)f''(x)
$$
for $x\in I^{\circ}$ and $f\in\mathcal{D}_{A^*}$.  Also
$0$ is an entrance boundary for $X^*$.  If $1$ is a reflecting boundary of $X$, then 1 is a regular absorbing boundary of $X^*$.  If $1$ is an entrance boundary of $X$, then 1 is an exit boundary of $X^*$.  
\ent
\begin{remark}
In all of the cases explored in Theorem \ref{thm:entrefl}, the domain of the dual generator is the same as the domain $\domas$ specified in Remark \ref{rem:gen}.
\end{remark}
\begx
\label{ex:bes}
\emph{ For $\alpha\geq0$, a
diffusion $X$ on $[0,1]$ is said to be a Bessel process with 
parameter $\alpha$ [written Bes($\alpha$)], reflected at $1$, if the generator of $X$ has the form $$A=\frac{1}{2} \frac{d^2}{dx^2}+\frac{\alpha-1}{2x} \frac{d}{dx},$$ and if for $f\in\doma$ we have $\frac{df}{dS}^-(1)=0$.   The behavior at the boundary $0$ is determined by the value of $\alpha$.  For $0<\alpha<2$, the boundary~$0$ is a regular reflecting boundary, and for $\alpha\geq 2$ the boundary~$0$ is an entrance boundary.  For our discussion of duality, we do not consider the case $\alpha=0$, for which~$0$ is an absorbing boundary.  
For $\alpha>0$, a simple application of Theorem \ref{thm:entrefl} gives that if $X$ is a Bes($\alpha$) process on $[0,1]$ with instantaneously reflecting behavior at $1$ begun in $\pi^{(x)}$, then $X^*$ is a Bes($\alpha$+2) process begun in $\delta_x$ absorbed at $1$.  In particular, the dual of reflecting Brownian motion, \ie,\ the Bes(1) process reflected at~$1$, is the Bes(3) process absorbed at~$1$.  For an extensive background treatment of Bessel processes, see \cite[Chapter4.3]{knight} or \cite[Chapter V--VI]{rogwil}.
}\enx

\begx
\emph{
For a second example, we turn to the 
Wright--Fisher gene frequency model from population genetics.  The Wright--Fisher diffusion $X$ is a diffusion on $[0,1]$ with generator of the form
$$A=\frac{1}{2} x(1-x) \frac{d^2}{dx^2}+[\alpha(1-x)-\beta x] \frac{d}{dx}.$$  The behavior at the boundaries is determined by the values of~$\alpha \geq 0$ and~$\beta \geq 0$.  We have that
 $$
0 \text{ is a(n) } 
\begin{cases}
\text{entrance boundary if } \alpha\geq 1/2,\\
\text{reflecting regular boundary if }0<\alpha<1/2,\\
\text{exit boundary if } \alpha=0,
\end{cases}
$$
and 
 $$
1 \text{ is a(n) } 
\begin{cases}
\text{entrance boundary if } \beta\geq 1/2,\\
\text{reflecting regular boundary if }0<\beta<1/2,\\
\text{exit boundary if } \beta=0.
\end{cases}
$$
If $X$ is a Wright--Fisher diffusion with $\alpha>0$ and $\beta=1/2$, then from \refT{thm:entrefl} we have that the strong stationary dual of $X$ is a Wright--Fisher diffusion with 
$\alpha^*=\alpha+(1/2)$ and $\beta^*=0$.  
For an extensive background on the Wright--Fisher model and its many applications, see \cite[Section 15.8]{kt2} or \cite[Chapter 10]{etku}.
}
\enx

Not surprisingly, we can also recover a partial converse to \refP{prop:sgl}.

\begl
Let $X$ and $X^*$ be diffusions on $[0, 1]$ and let~$0$ and $1$ be either instantaneously reflecting or entrance boundaries for $X$.  Then an intertwining 
\begin{equation}\label{eq:link2}
\gL T_t=T_t^*\gL\ \mathrm{(}\text{for all }t\geq0\mathrm{)}
\end{equation}
of the one-parameter semigroups by the link~$\gL$
together with the initial condition~\eqref{eq:init} implies that $X^*$ is a strong stationary dual of $X$.
\enl

\begin{proof}
Let $f\in\doma$.  Then by the infinitismal characterization of the generator $A$ in terms of its associated semigroup $(T_t)_{t=0}^{\infty},$ we have that 
$$Af=\mathrm{lim}_{t\downarrow 0}\frac{T_tf-f}{t},$$
with convergence in the uniform norm.
Now
\begin{align*}
(\gL A f)(x)
&= \int_0^x\!\pi^{(x)}(y)\left[ \lim_{t\downarrow 0}\frac{T_tf(y)-f(y)}{t} \right]\,dy\\
&= \lim_{t\downarrow 0}\int_0^x\!\pi^{(x)}(y)\,\frac{T_tf(y)-f(y)}{t}\,dy\\
&= \lim_{t\downarrow 0}\frac{(\gL T_t f)(x) - (\gL f)(x)}{t}\\
&=\lim_{t\downarrow 0}\frac{(T_t^* \gL f)(x) - (\gL f)(x)}{t}\\
&= (A^* \gL f)(x),
\end{align*}
where the last limit's existence is guaranteed by that of the first.  
As $\gL Af\in C(I)$, it follows that $A^*\gL f \in C(I)$.  This implies that 
$$A^*\gL f=\mathrm{lim}_{t\downarrow 0}\frac{T^*_t\gL f-\gL f}{t},$$
in the uniform norm
(see \cite[Theorem 7.7.3]{lamp}), and $\gL f\in \domas.$ 
Combined, we have that $\gL|_{\mathcal{D}_A}\subset\domas$, and that on $\doma$ we have $\gL A=A^* \gL$ as desired. 
\end{proof}

\begr
Intertwinings of Markov semigroups
have been well studied, appearing for example in \cite{dynk}, \cite{ROGPIT}, etc.  
In the context of (\ref{eq:link2}), the transition operator $\Lambda$ 
is the following Markov kernel from $[0,1]$ to $[0,1]$: For $x\in[0,1]$ and $A\in\mathcal{B}$ we have  
$$
\Lambda(x,A) =  \Pi^{(x)}(A).
$$
\enr

\begr
\label{rem:inter}   
If~\eqref{eq:link2} holds, then
$$
(U_t\Pi_0,f) =(\Pi_0,T_tf) = (\Pi_0^*,\gL T_tf) = (\Pi^*_0,T_t^*\gL f) = (U_t^*\Pi_0^*,\gL f),
$$
mirroring 
the corresponding result that algebraic duality via link $L$ of Markov chains yields $\pi_t = \pi^*_t L$.
\enr

\section{Approximating duality via Markov chains}
\label{sec:2}
The purpose of the present section is twofold.  Presently suppressing all details (which will be spelled out in full detail later in the section), we will show that a suitably defined sequence of Markov chains $X^{\gD}$ and their corresponding strong stationary duals $\widehat{X}^{\gD}$, as defined in \cite{DFSST}, converge respectively to our primal diffusion $Y=S(X)$ (in natural scale) and its strong stationary dual $Y^*$.  By establishing the newly defined diffusion strong stationary dual as a limit of an appropriately defined sequence of classical Markov chain strong stationary duals, we ground our definition and our present work in the classical theory. 

In addition to tethering our duality to the classical theory, this has a number of interesting consequences.  For example, we believe one of the great triumphs of strong stationary duality was its application in the perfect sampling algorithms of \cite{fps} and \cite{MR2001m:60164}.  Via the work in the present section, for our primal diffusion $Y$ we could approximately sample perfectly from $\Pi_Y$ by using the theory of \cite{fps} to perfectly sample from the stationary distributions of the approximating sequence of chains.  We could also use our approximating sequence of chains to study cut-off type behaviors of the dual hitting times of state $S(1)$, and hence of the primal diffusion's separation distance from stationarity.  We are also able to recover the 
dual-hitting-time/primal-separation duality of the classical Markov chain theory in the diffusion setting by passing to appropriate limits; see Section \ref{S:SST} for full details.

This section is laid out as follows:  First assuming instantaneously reflecting boundaries for our primal diffusion $Y$, in Sections \ref{sec:pc}--\ref{sec:dualconv} we explicitly spell out the one-dimensional convergence of our primal and dual sequences of Markov chains to the corresponding primal and dual diffusions.  In Section \ref{sec:ent}, we prove the corresponding convergence theorems in the case when our primal diffusion has entrance boundaries at $0$ and/or~$1$.

\subsection{Primal convergence}
\label{sec:pc}
Let $D_{I}[0,\infty)$ be 
the space of cadlag functions from $[0,\infty)$ into $I$.  We can equip  $D_{I}[0,\infty)$ with a metric $d$ defined by
$$d(x,y)=\inf_{\gl\in B}\bigg[\left(\sup_{s>t\geq 0}\bigg|\log\frac{\gl(s)-\gl(t)}{s-t}\bigg|\right) \vee \int_0^{\infty}e^{-u}d(x,y,\gl,u)\,du\bigg]$$
where $B$ is the set of strictly increasing Lipschitz continuous functions from $[0,\infty)$ to $[0,\infty)$ with the additional property that 
$$
\gl \in B\mbox{\ implies\ }\sup_{s>t\geq 0}\bigg|\log\frac{\gl(s)-\gl(t)}{s-t}\bigg|<\infty,
$$
and $$d(x,y,\gl,u):=\sup_{t\geq 0}\left(|x(t\wedge u)-y(\gl(t)\wedge u)|\wedge 1\right).$$
The topology induced by $d$ is known as the \textit{Skorohod topology}, and under this topology $D_{I}[0,\infty)$ is both complete and separable (as $I$ is both complete and separable).  For more background on $D_{I}[0,\infty),$ see \cite[Sections 3.5--3.10]{etku} 
or~\cite[Chapters 2--3]{billing}.  

We will consider stochastic processes with sample paths in $D_{I}[0,\infty)$ as $D_{I}[0,\infty)$-valued random variables and we will say that $X_n\Rightarrow X$ if we have convergence in law of the corresponding $D_{I}[0,\infty)$-valued random variables. Note that $X_n \Rightarrow X$ implies convergence of the associated 
finite-dimensional distributions of $X_n$ to those  of~$X$ (see \cite[Theorem 3.7.8]{etku}), \ie,\ for all $\{t_1,\ldots,t_m\} \subset \{t\geq 0\,|\,\p(X(t)=X(t-))=1\}$ we have
$$(X_n(t_1),\ldots,X_n(t_m))\Rightarrow(X(t_1),\ldots,X(t_m)).$$  

As in \refS{S:Section 1}, let $X$ be a regular diffusion on $I$ with instantaneous reflection at the boundaries of 
$I$ and  scale function $S = S_X$.  To ease exposition, we will consider $Y=S_X(X)$, the regular diffusion in natural scale on $\scrs = [S_X(0), S_X(1)]$ (note that $0$ and $1$ being instantaneously reflecting implies that $S_X(0)$ and $S_X(1)$ are finite), and assume $S_Y$ has been scaled to make $s_Y\equiv 1$.  
The speed function of $Y$ is $M_Y=M_X\circ S_X^{-1}:\scrs^{\circ}\rightarrow \mathbb{R}$ (where $M_X$ is the speed function of $X$).  As with $X$, define the speed measure of $Y$ as the nonnegative measure on $\scrs^{\circ}$, denoted $M_Y(\cdot)$, defined via $M_Y(x,y]=M_Y(y)-M_Y(x)$.

It follows easily (from the analogous results for $X$) that $N_Y(S(0))<\infty$ and 
$\Sigma_Y(S(0))<\infty$ [with analogous results holding at $S(1)$], and therefore the boundaries of $\scrs$ are regular for $Y$.  Regularity of $X$ and the fact that $S$ is continuous and strictly increasing on $I$ implies that $Y$ is regular and that the boundaries of $\scrs$ are instantaneously reflecting for $Y$.

The generator of~$Y$ can be expressed as $(A_Y f)(y) = \frac{1}{2} b_Y (y) f''(y)$ with 
$b_Y(y) = b_X(x) s_X^2(x)$ where $y = S_X(x)$. Note 
that $M_Y(\scrs^{\circ})=M_X(I^{\circ})<\infty$ and so there exists a unique invariant measure for $Y$ which we will denote $\Pi_Y$.  Observe 
\begin{align*}
M_Y((c,d])&=M_Y(d)-M_Y(c)=M_X(S^{-1}(d))-M_X(S^{-1}(c))\\
&=\int_{S^{-1}(c)}^{S^{-1}(d)}\!m_X(z)\,dz=\int_c^d \frac{m_X(S_X^{-1}(w))}{s_X(S_X^{-1}(w))}\,dw
\end{align*} 
so that $M_Y$ (resp.,\ $\Pi_Y$) has density $m_Y(y) = m_X(S_X^{-1}(y)) / s_X(S_X^{-1}(y))$
(resp., density $\pi_Y=\alpha\,m_Y$ for some constant $\alpha$).  
On $\scrs^{\circ}$, $m_Y=b_Y^{-1}$ implies that $\pi_Y b_Y$ is constant.  Assume that $b_Y$ can be extended to a function in $C(\scrs)$, so that $b_Y(S(0))$ and $b_Y(S(1))$ are well defined, and similarly extend $\pi_Y b_Y$ to $C(\scrs)$ via
$\pi_Y(z)b_Y(z)=\alpha$ for $z \in \{S(0), S(1)\}$.

For the remainder of the section, we shall be working with the diffusion $Y$ rather than $X$, and so 
we will drop the $Y$ subscript from $b_Y$, $\pi_Y$, $\Pi_Y$, etc.

For $n \in \{2,3,4,\ldots\}$, define $\Delta=\Delta_n=[S(1)-S(0)]/n$.  
As in \cite[Chapter~6]{bhat}, define a birth-and-death transition matrix $P^{n}$ on state space
$$\scrs^{n} := \{S(0),\  S(0)+\Delta,\ S(0)+2\Delta,\ldots,S(1)-\Delta,\ S(1)\}$$
by setting (for ease of notation, we write~$i$ for $S(0)+i\Delta$ here):
$$
P^{n}(i,i+1)=P^{n}(i,i-1):=\frac{b(i)h}{2\Delta^2} \quad \text{for }0<i<n\text{ and}
$$
$$
P^{n}(0,1):=\frac{b(0)h}{\Delta^2},\ \ P^{n}(n,n-1):=\frac{b(n)h}{\Delta^2};
$$
here 
\begin{equation}
\label{eq:h}
h = h_{n} := \frac{\Delta^2}{2 \sup_y b(y)}
\end{equation}
 is chosen to make $P^{n}$ monotone.  

Note that for $i\in\{1,\ldots,n-1\}$ we have 
$$\pi(i)P^{n}(i,i+1)=\pi(i+1)P^{n}(i+1,i),$$
and at the boundaries we have
$$\pi(0)P^{n}(0,1)=2\pi(1)P^{n}(1,0),\ \pi(n)P^{n}(n,n-1)=2\pi(n-1)P^{n}(n-1,n).$$
It follows that there exists a constant $C^{n}$ such that 
\begin{equation}
\label{eq:statprim}
\pi^{n}(i)=\begin{cases}C^{n}\pi(i), & i=1,\ldots,n-1;\\
C^{n} \pi(i) / 2, & i=0,\ n
\end{cases}
\end{equation}
is the unique invariant probability distribution for $P^{n}$.  

Let $\pi_0^{n}$ be a probability measure on $\scrs^{n}$, and let $P^{n}$ be the transition matrix for a discrete-time birth-and-death chain $X^{n}$, begun in $\pi_0^{n}$, on state space $\scrs^{n}$ [we write $X^{n}\sim (\pi_0^{n},P^{n})$ as shorthand].

\begt
\label{thm:mccon}
Assume there exists a constant $\gd>0$ such that $b \geq \gd$ everywhere and that we can continuously extend $b$ to the boundaries of $\scrs$.
For $n=2,3,4,\ldots,$ define the continuous-time stochastic process $Y^{n}$ by setting 
$Y^{n}_t := X^{n}_{\lfloor t/h_n\rfloor}$ for $t \geq 0$.  If $Y^{n}_0 \Rightarrow Y_0$, 
then 
$Y^{n} \Rightarrow Y$. 
\ent

Our main proof tool will be the following theorem, adapted from \cite[Corollary 4.8.9 and Theorem 1.6.5]{etku}:
\begin{theorem}
\label{thm:converge}
Let~$A$ be the generator (as in Section \ref{S:Section 1}) of a regular diffusion process~$Y$ with state space~$\mathcal{Y}$.  Assume $h_n > 0$ converges to~$0$ as $n \to \infty$.  Let $X^{n}\sim(\pi_0^{n},P^{n})$ be a Markov chain on metric state space $\mathcal{Y}^{n} \subset \mathcal{Y}$ and define $Y^{n}_t := X^{n}_{\lfloor t/h_{n} \rfloor}$.  Further assume $Y^{n}_0\Rightarrow Y_0$.  Letting $\mathcal{B}(\scry^n)$ be the space of real-valued bounded measurable functions on $\scry^n$, define $T^{n}:\mathcal{B}(\scry^{n})\rightarrow \mathcal{B}(\scry^n)$ via
$$
T^{n} f(x)=\e_x f(X^{n}_1).
$$ 
Let $\rho_{n}:C(\scry)\rightarrow \scrb(\scry^{n})$ be defined via $\rho_{n}f(\cdot)=f |_{\scry^{n}}(\cdot)$.  If $\doma$ is an algebra that strongly separates points, and
\begin{equation}
\label{eq:corem}
\lim_{n\rightarrow \infty} \sup_{y\in\scry^{n} }\left| (A^{n}\rho_{n} f)(y) - (A f)(y) \right|=0
\end{equation}
for 
all $f\in \doma$, then $Y^{\gD}\Rightarrow Y$.
\end{theorem}

The adaptation of Theorem \ref{thm:converge} from \cite[Corollary 4.8.9 and Theorem 1.6.5]{etku} is spelled out explicitly in Appendix~\ref{section:extpf}, as the notation between \cite{etku} and the present section differs considerably.
\begin{proof}[Proof of \refT{thm:mccon}]
Let $f\in\doma$, so that $(A f)(y)=\frac{1}{2}b(y)f''(y)$.  Using 
$$
\doma=\{f\in C(\scrs)\cap C^2(\scrs^{\circ})\ |\ Af\in C(\scrs),\ f'(S(0)+)=f'(S(1)-)=0\},
$$ 
we find that
$b f'' \in C(\scrs)$.  As $b(y) \geq \delta>0$ for all $y\in\scrs$, we have 
$1/b\in C(\scrs)$ and therefore
$f''\in C(\scrs)$.  This implies that 
$$
\doma=\{f\in C^2(\scrs)\ |\ Af\in C(\scrs),\ f'(S(0)+)=f'(S(1)-)=0\},
$$  
which is indeed an algebra that strongly separates points.
It follows that 
\begin{align*}
&\lim_{n\rightarrow \infty} \sup_{y\in\scrs^{n}\setminus \{S(0), S(1)\} }\left|h^{-1} ((T^{n} - I) f)(y)- \frac{1}{2}b(y) f''(y)\right|=\\
&\lim_{n\rightarrow \infty} \sup_{y\in\scrs^{n}\setminus \{S(0), S(1)\} }\left|\frac{1}{2}b(y)\frac{f(y+\Delta)-2f(y)+f(y-\Delta)}{\Delta^2}
-\frac{1}{2}b(y) f''(y)\right|=0.
\end{align*}
Likewise, 
\begin{align*}
h^{-1} ((T^{n} - I) f)(S(0))&=b(S(0))\frac{f(S(0)+\Delta)-f(S(0))}{\Delta^2}
\to \frac{1}{2}b(S(0))f''(S(0)), \\
h^{-1} ((T^{n} - I) f)(S(1))&=b(S(1))\frac{f(S(1)-\Delta)-f(S(1))}{\Delta^2}
\to \frac{1}{2}b(S(1))f''(S(1)).
\end{align*}
Therefore 
$$
\sup_{y\in\scrs^{n} }\left| (A^{n} \rho_{n}f)(y) - (A f)(y) \right| \to 0,
$$ 
establishing~\eqref{eq:corem}, and the result follows.
\end{proof}

\begr
\label{remark:initconv}
Recall that $\Delta=\Delta_n=[S(1)-S(0)]/n$.
For $x\in(S(0),S(1))$, let $i_{n,x} := \left\lfloor [x-S(0)] / \gD \right\rfloor$, and denote the invariant measure $\pi^{n}$ truncated to $\{S(0),S(0)+\gD,\ldots,S(0)+i_{n,x}\gD\}$ by $\pi^{n, i_{n,x}}$.  If $Y^n_0\sim\pi^{n, i_{n,x}}$ and $Y$ is begun with density 
$\pi^{(x)}$, then for $y\in(S(0)+k\gD,S(0)+(k + 1)\gD)$ with $0 \leq k <i_{n,x}$ 
we have
$$
\p^{n}(Y^{n}_0\leq y)=\frac{\gD\sum_{j=0}^k \pi^{n}(S(0)+j\gD)}{\gD\sum_{j=0}^{i_{n,x}} \pi^{n}(S(0)+j\gD)}\rightarrow\frac{\int_{S(0)}^y\!\pi(z)\,dz}{\int_{S(0)}^x\!\pi(z)\,dz} = \p(Y_0 \leq y).
$$
If $X^{n}$ is begun in $\pi^{n,i_{n,x}}$, it follows that $X^{n}_0=Y^{n}_0\Rightarrow Y_0$.

If instead~$Y$ is begun deterministically at $S(0)$, then letting 
$X^{n}_0=Y_0^n=S(0)$ trivially yields $X^{n}_0=Y^{n}_0\Rightarrow Y_0$.
\enr

\begr
\label{remark:conddual} For the sequence of birth-and-death chains $(\pi_0^{n},P^{n})$, where either 
$\pi_0^{n}=\delta_{S(0)}$ for 
each~$n$ or $x \in (S(0), S(1)]$ is given and fixed and $\pi_0^{n}=\pi^{n, i_{n,x}}$ for each $n$, we are guaranteed the existence of a sequence of birth-and-death strong stationary dual chains by \cite[eqs.\ (4.16a)--(4.16b)]{DFSST} because of the following two observations.

(a)~$P^{n}$ is 
monotone.  Indeed, for $i = 0, \dots, n - 1$ (again employing the shorthand $i$ for 
$S(0)+i\gD$) we easily see
$$
P^{n}(i, i + 1) + P^{n}(i + 1, i) \leq 1.
$$ 

(b)~The ratio $\pi^{n}_0 / \pi^{n}$ of initial probability mass function to the 
stationary distribution
is non-increasing. 
\enr

\subsection{Dual convergence}
\label{sec:dualconv}
As in \cite{DFSST}, construct on the same probability space  
as for $X^{n}$ a strong stationary dual 
$\widehat{X}^{n}\sim (\hat{\pi}_0^{n},\widehat{P}^{n}$) of $X^n$ using the 
link~$\Lambda^n$ of truncated stationary distributions 
(here and below, for ease of notation, $i$ is again used as shorthand for 
$S(0)+i\gD$): 
$$
\gL^{n}(i,j) := \mathbbm{1}\{j\leq i\}\frac{\pi^{n}(j)}{H^{n}(i)};
$$  
we have used the shorthand $H^{n}(i):=\sum_{j=0}^i \pi^{n}(j)$.  Note that $\widehat{X}^{n}$ is also a birth-and-death chain on $\scrs^{n}$.

We further assume that for all $n$,
there is a fixed $x\in(S(0),S(1)]$ such that
$X^{n}_0\sim \pi^{n, i_{n,x}}$ (so that $\widehat{X}^{n}_0=i_{n,x}$).  The one-step transition matrix $\Ph^{n}$ for  
$\widehat{X}^{n}$ is given by 
\begin{equation}
\label{eq:death}
\widehat{P}^{n}(i,i-1)=\frac{H^{n}(i-1)}{\HD(i)}\frac{b(i)h}{2\gD^2}=\frac{b(i)h}{2\gD^2}-\frac{h\cdot \alpha\cdot\CD}{\HD(i)2\gD^2}\ \ \ \text{ for } 0<i<n,
\end{equation}
\begin{equation}
\label{eq:birth}
 \widehat{P}^{n}(i,i+1)=\frac{H^{n}(i+1)}{\HD(i)}\frac{b(i+1)h}{2\gD^2}=\frac{h\cdot \alpha\cdot\CD}{\HD(i)2\gD^2}+\frac{b(i+1)h}{2\gD^2}\ \ \ \text{ for }0<i < n,
\end{equation}
\begin{equation}
\label{eq:from0}
\widehat{P}^{n}(0,1)=\frac{H^{n}(1)}{\HD(0)}\frac{b(1)h}{\gD^2},
\end{equation}
\begin{equation}
\label{eq:absorb}
\widehat{P}^{n}(n,n)=1,
\end{equation}
with $\widehat{P}^{n}(i,i)$ for $0 \leq i < n$ defined so that the rows of~$\widehat{P}^{n}$ sum to unity. 
We 
next show the following theorem:
\begt
\label{thm:dualconv}
With assumptions as in Theorem \ref{thm:converge}, further
assume that $b\in C^2(S(0),S(1)]$, and $\pi(S(0))>0$ and $\pi(S(1))>0$.  
  Define the continuous-time processes $(\widehat{Y}^{n}_t) := (\widehat{X}^{n}_{\lfloor t/h_n \rfloor})$ [with $h_n$ defined as at (\ref{eq:h})], and assume $\widehat{Y}^{n}_0\Rightarrow Y^*_0$.  Then 
$$\widehat{Y}^{n}\Rightarrow Y^*,$$
where $Y^*$ is a SSD of~$Y$ in the sense of Definition~\ref{D:1.1}.  
\ent

We will prove \refT{thm:dualconv} after a series of preliminary results.  We begin by putting $Y^*$ into natural scale, \ie, we consider the diffusion $Z^*=S^*(Y^*)= -1 / M(Y^*)$ on state space $\scrs^*=(S^*(S(0)),S^*(S(1))]=(-\infty,S^*(S(1))]$, as $S(0)$ is an entrance boundary for the dual diffusion $Y^*$.
Note that the infinitesimal parameters of $Z^*$ are given on $(-\infty,S^*(S(1)))$ [recalling~\eqref{eq:dualscale}] by
$$a_{Z^*}\equiv 0,\ \ b_{Z^*}(S^*(y))=b(y)s^*(y)^2=\frac{b(y)m^2(y)}{M^4(y)}=\alpha^2\frac{b(y)\pi^2(y)}{\Pi^4(y)}$$  (recall $\alpha$ is the constant such that $\pi=\alpha\cdot m=\alpha/b$).  
Note also that under
the assumptions of Theorem~\ref{thm:dualconv}, we have $b\in C^2(S(0),S(1)]$ (and so 
$\pi(\cdot)\propto b^{-1}(\cdot)$ implies that $\pi\in C^2(S(0),S(1)$] as well), and also $\pi>0$ on $\partial\scrs$.  Note also that
$$
(S^*)'(i)=\alpha \frac{\pi(i)}{\Pi^2(i)},
$$
$$
(S^*)''(i)=\alpha\frac{\Pi^2(i)\pi'(i)-2\pi^2(i)\Pi(i)}{\Pi^4(i)},
$$
\begin{equation*}
\label{eq:c3}
(S^*)'''(i)=\alpha\frac{\Pi^2(i)\pi''(i)-2\pi'(i)\Pi(i)\pi(i)}{\Pi^4(i)}-\alpha\frac{4\Pi^3(i)\pi(i)\pi'(i)-6\pi^3(i)\Pi^2(i)}{\Pi^6(i)},
\end{equation*}
which implies that $S^*\in C^3[-S(0)+i_0\Delta,S(1)]$ for any $i_0>0$. 

Define
$$\widehat{Z}_t^{{n}}:=S^*\left(\widehat{X}_t^{{n}}\right),$$
and note that this is a birth and death chain on state space $$\scrs^{*,n}:=\{S^*(S(0)),\ S^*(S(0)+\gD),\ldots,S^*(S(1)-\gD),S^*(S(1))\}.$$  
For all~$x$ of the form $S^*(S(0)+i\gD)$, where $i > 0$ is an integer, define
\begin{equation}
\label{eq:driftdel}
a_{\widehat{Z}^{n}}(x):=\frac{1}{h_{n}}\e_x\left(\widehat Z_1^{n}-x\right).
\end{equation}
Then we have the following proposition.
\begp
\label{prop:driftcon}
With 
the  assumptions of Theorem \ref{thm:dualconv} and assuming $i>0$, letting $R<\infty$ be fixed, it follows 
that
$$
\lim_{n \to \infty}\sup_{x\in\scrs^{*,n}: |x|<R}|a_{\widehat{Z}^{n}}(x)|=0.
$$ 
\enp
\begin{proof} As before, we abbreviate $S(0)+i\gD$ as $i$, so that $x$ can be written $x=S^*(i)$.  Let
$${\rm A}:=\frac{\bss(i+1)-2\bss(i)+\bss(i-1)}{\gD^2}\frac{b(i)}{2};\ \ {\rm B}:=\frac{\bss(i+1)-\bss(i)}{\gD}\frac{b(i+1)-b(i)}{2\gD};$$
$${\rm C}:=\frac{\alpha\cdot C^{n}}{\gD \HD(i)}\frac{\bss(i+1)-\bss(i-1)}{2\gD}.$$
Then \eqref{eq:death}--\eqref{eq:birth}
allow us to rewrite (\ref{eq:driftdel}) as 
$$a_{\widehat{Z}^{n}}(x)=\rm{A}+B+C$$
for $x\neq S^*(S(1))$.

There exist constants $\gamma$ and~$\gd$ such that for all $|x|=\bss(i)<R$ we have 
$$0<\gd<\inf_n i\gD<\sup_n i\gD<\gamma<\infty.$$  A Taylor expansion of $\bss(\cdot)$ combined with 
$\bss\in C^3[\gd,S(1)]$ gives
$$\lim_{n \to \infty}\,\sup_{x\in\scrs^{*,n}:|x|<R,\,x\neq S^*(S(1))}\bigg|{\rm A}-\alpha\frac{\Pi^2(i)\pi'(i)-2\pi^2(i)\Pi(i)}{\Pi^4(i)}\frac{b(i)}{2}\bigg|=0,$$
or equivalently
$$\lim_{n \to \infty}\,\sup_{x\in\scrs^{*,n}:|x|<R,\,x\neq S^*(S(1))}\bigg|{\rm A}-\frac{\alpha^2}{2}\left(\frac{\pi'(i)}{\pi(i)\Pi^2(i)}-\frac{2\pi(i)}{\Pi^3(i)}\right)\bigg|=0;$$
and 
$$\lim_{n \to \infty}\,\sup_{x\in\scrs^{*,n}:|x|<R,\,x\neq S^*(S(1))}\bigg|{\rm B}-\alpha \frac{\pi(i)}{\Pi^2(i)}\frac{b'(i)}{2}\bigg|=0,$$
or equivalently
$$\lim_{n \to \infty}\,\sup_{x\in\scrs^{*,n}:|x|<R,\,x\neq S^*(S(1))}\bigg|{\rm B}+\alpha^2 \frac{\pi'(i)}{2\Pi^2(i)\pi(i)}\bigg|= 0.$$

The function $\pi$ is uniformly continuous on $\scrs$, and, since the Riemann sum of a uniformly continuous function converges uniformly to the corresponding Riemann integral, we 
have for the sup-norm $\| \cdot \|_{\infty}$ on $\scrs$ that
$$
\lim_{n \to \infty}||(C^{n})^{-1}\gD \HD-\Pi||_{\infty}= 0.
$$  
By regularity of the primal diffusion $Y$, we have $\Pi(i)>0$ for $i\gD>\delta>0$, and therefore for 
such~$i$ we have $(C^{n})^{-1}\gD \HD(i)>0$ and also $(C^{n})^{-1}\gD \HD(i) \Pi(i)>0$.  Note that $(C^{n})^{-1}\gD \HD(\cdot)\Pi(\cdot)$ is a bounded increasing function in $i$.  
All of this leads to 
\begin{align*}
&\lim_{n \to \infty}\,\sup_{x\in\scrs^{*,n}:|x|<R,\,x\neq S^*(S(1))}\bigg|\frac{C^{n}}{\gD \HD(i)}-\frac{1}{\Pi(i)}\bigg|=\\ 
&\lim_{n \to \infty}\,\sup_{x\in\scrs^{*,n}:|x|<R,\,x\neq S^*(S(1))}\bigg|\frac{C^{n}}{\gD \HD(i)\Pi(i)}\bigg|\cdot\bigg|\Pi(i)-\frac{\gD\HD(i)}{C^{n}}\bigg|=0.\end{align*}
Therefore
$$\lim_{n \to \infty}\,\sup_{x\in\scrs^{*,n}:|x|<R,\,x\neq S^*(S(1))}\bigg| {\rm C}-\frac{\alpha^2\pi(i)}{\Pi(i)^3}\bigg|= 0.$$
Combining our results for A, B, and C with the observations that 
$a_{\widehat{Z}^n}(S^*(S(1)))=0$, we find
$$
\lim_{n \to \infty}\,\sup_{|x|<R}|a_{\widehat{Z}^{n}}(x)|= 0,
$$ 
as desired.
\end{proof}

Next, define
\begin{equation}
\label{eq:diffdel}
b_{\widehat{Z}^{n}}(x):=\frac{1}{h_{n}}\e_x \left( \widehat Z^{n}_1 -x \right)^2,
\end{equation}
where again we require 
that $x=S^*(S(0)+i\gD)$ for some positive integer $i$.

\begp
\label{prop:diffcon}
With the assumptions of Theorem \ref{thm:dualconv}, letting $R<\infty$ be fixed, we have
$$\lim_{n \to \infty}\,\sup_{x\in\scrs^{*,n}:|x|<R}\bigg|b_{\widehat{Z}^{n}}(x)-b_{Z^*}(x)\bigg|= 0,$$
where $b_{Z^*}(S^*(S(1)))=0$ by the absorbing behavior of the boundary at $S^*(S(1))$.
\enp
\begin{proof}
We have
$$b_{Z^*}(x)=\alpha^2\frac{\pi^2((\bss)^{-1}(x))}{\Pi^4((\bss)^{-1}(x))}b((\bss)^{-1}(x))$$ for $x\neq S^*(S(1))$.
There exists a constant $\gd > 0$ such that for all $\gD$ and all~$x$ satisfying $|x| < R$, if we write $x=S^*(S(0)+i \gD)$ 
then $i \gD \geq \gd$.
Let 
$$ {\rm A}:=(\bss(i+1)-\bss(i))^2\frac{\alpha\cdot C^{n}}{2\gD^2\HD(i)};\ \ {\rm B}:=\frac{(\bss(i+1)-\bss(i))^2}{\gD^2}\frac{b(i+1)}{2} ;$$
and
$${\rm C}:=\frac{(\bss(i-1)-\bss(i))^2}{\gD^2}\frac{b(i)}{2};\ \ {\rm D}:=-(\bss(i-1)-\bss(i))^2\frac{\alpha\cdot C^{n}}{2\gD^2\HD(i)}.$$
Note that $b_{\widehat{Z}^{n}}(x)=\rm A+B+C+D$ for $x\neq S^*(S(1))$.  

As 
in the proof of Proposition~\ref{prop:driftcon},  
$$\lim_{n \to \infty}\,\sup_{x\in\scrs^{*,n}:|x|<R,\,x\neq S^*(S(1))}\bigg| {\rm B}-\alpha^2\frac{\pi^2(i)}{2\Pi^4(i)}b(i)\bigg|= 0;$$
$$\lim_{n \to \infty}\,\sup_{x\in\scrs^{*,n}:|x|<R,\,x\neq S^*(S(1))}\bigg| {\rm C}-\alpha^2\frac{\pi^2(i)}{2\Pi^4(i)}b(i)\bigg|=0.$$
Rewrite A (with analogous results holding for D) as 
$${\rm A}=(\bss(i+1)-\bss(i))\frac{\bss(i+1)-\bss(i)}{\gD}\frac{\alpha\cdot C^{n}}{2\gD\HD(i)}.$$
From the uniform continuity of $(\bss)''(\cdot)$ on bounded intervals, $\frac{\bss(i+1)-\bss(i)}{\gD}$ converges uniformly for $|x|=|\bss(i)|<R$ to $(\bss)'(i)$, which is uniformly bounded for $|x|=|\bss(i)|<R$.  Also, 
$(C^{n})^{-1}2\gD\HD(i)$ is bounded away from 0 for $|x|=|\bss(i)|<R$, and $\bss(\cdot)$ is uniformly continuous for $|x|=|\bss(i)|<R$.  Hence
$$\lim_{n \to \infty}\,\sup_{x\in\scrs^{*,n}:|x|<R,\,x\neq S^*(S(1))}|{\rm A}|= 0; \qquad 
\lim_{n \to \infty}\,\sup_{x\in\scrs^{*,n}:|x|<R,\,x\neq S^*(S(1))}| {\rm D}|=0.$$
Lastly, note that $b_{\widehat{Z}^{n}}(S^*(S(1)))=0$,
which finishes the proof.
\end{proof}


We are now ready to prove Theorem \ref{thm:dualconv}:
\begin{proof}[Proof of \refT{thm:dualconv}]


With $x>S(0)$ fixed so that $\p(Z_0^*)^{-1}=\delta_{\bss(x)}$,
let $R$ be such that $|\bss(x)|<R$ and $|\bss(S(1))|<R$, and define
$$\tau_R:=\inf\{t\geq 0\ :\  |Z^*_t|=R\}.$$
It follows that $Z^*(\cdot\wedge\tau_R)$ is equal in distribution to the diffusion process with state space 
$\scrs^*_R:=[-R,S^*(S(1))]$ and generator
$$A^*_R:=\frac{1}{2} b_{Z^*}(\cdot)\frac{\partial^2}{\partial x^2}$$ operating on the domain
$$\mathcal{D}_R:=\left\{f\in C(\scrs^*_R)\cap C^2[(\scrs^*_R)^{\circ}]\,\big|\,A^*_R f\in C(\scrs^*_R),\ A^*_R f(-R)=A^*_R f(\bss(S(1)))=0 \right\}.$$
Since $b_{Z^*}>0$ on $\scrs^*$, it follows that if $f\in\mathcal{D}_R$ then $f\in C^2(\scrs^*_R)$.  Let 
$$\tau_R^{n}:=\inf\{t\geq 0\ :\  |\widehat{Z}^{n}_t|\geq R\},$$
and define the sequence of absorbing Markov chains
$\widehat{V}^{n}(\cdot):=\widehat{Z}^{n}(\cdot\wedge\tau_R^{n})$ with state space
$$\scrs^{*,n}_R:=\{x\in \scrs^{*,n}: x \geq \lceil R \rceil_\gD\}$$  
where $\lceil y \rceil_{\gD}$ ``rounds'' $y$ to the smallest element $\geq y$ in the grid 
$\{S^*(S(0)),\,S^*(S(0)+\gD),\ldots, S^*(S(1)-\gD),\,S^*(S(1))\}$.
Lastly define $V^{*,n}(t):=\widehat{V}^{n}(\lfloor t/h\rfloor)$.

From \cite[Corollary 4.8.9]{etku}, to prove that $V^{*,n}$ converges (as $n \to \infty$) to $Z^*(\cdot\wedge\tau_R)$, it suffices 
to show that for each fixed $f\in \mathcal{D}_R$ we have
\begin{equation}
\label{3.5A}
\lim_{n \to \infty}\sup_{x\in \scrs^{*,n}_R}|\rho_{n}f(x)-f(x)| = 0
= \lim_{n \to \infty}\sup_{x\in \scrs^{*,n}_R}|\rho_{n} A^*_R f(x)-\widehat{A}^{n}f(x)|
\end{equation}
where $\widehat{A}^{n}f(x) := h_{n}^{-1}[\e_x f(\widehat{V}^{n}_1) - f(x)]$.
The first equality in~\eqref{3.5A} is trivial.  Consider the second equality.  
At $x=-R$ (if $-R\in\scrs^{*,n})$ or $x=\bss(S(1))$, we 
have
$$
|\rho_{n}  A^*_R f(x)-\widehat{A}^{n}f(x)| = 0,
$$ 
since both $\rho_{n}  A^*_R f(x)$ and 
$\widehat{A}^{n}f(x)$ equal 0 for $x=-R$ (if $-R\in\scrs^{*,n})$ or $x=S^*(S(1))$.

For~$x$ in the interior of 
$\scrs^{*,n}_R$, from a Taylor expansion of~$f$ with remainder in intermediate-point form we find 
\begin{equation}
\label{3.51}
\left| \widehat{A}^{n}f(x) - 
\left[ f'(x) a_{\widehat Z^{n}}(x) + \frac{f''(x)}{2} b_{\widehat Z^{n}}(x) \right] \right| 
\leq \frac{c}{2} b_{\widehat Z^{n}}(x),
\end{equation}  
where, with $x = S^*(S(0) + i \gD)$, we take 
$$
c = \max_{\xi\in[S^*(i-1),S^*(i+1)]} |f''(\xi) - f''(x)|.
$$
From \eqref{3.51}, \refP{prop:driftcon}, \refP{prop:diffcon}, and the fact that $f\in C^2(\scrs^*_R)$, we have that 
$\widehat{A}^{n}f(x)$ converges uniformly to $A^*_Rf(x)$, and so \eqref{3.5A} is proven.

We have now established that  $V^{*,n}$ converges in distribution to $Z^*(\cdot\wedge\tau_R)$.
From \cite[Theorem 11.1.1]{SV}, we have that $\widehat Z^{n}\Rightarrow Z^*$.  Lastly, noting that $(S^*)^{-1}(\cdot)$ is well-defined and measurable (indeed it is continuous!) we have that
$$\widehat{Y}^{n}=(S^*)^{-1}(\widehat{Z}^{n})\Rightarrow (S^*)^{-1}(Z^*)=Y^*$$ by \cite[Theorem 3.10.2]{etku}.
\end{proof}
%

\subsection{Convergence extended to entrance boundary cases}
\label{sec:ent}
For $0$ an entrance boundary of $X$ and $1$ reflecting, again consider $Y=S(X)$, a regular diffusion in natural scale on $\scrs=[-\infty,S(1)]$ begun in $\Pi_0=\pi^{(x)}$ for fixed $x\in\scrs\setminus\{-\infty\}$, the stationary measure for $Y$ truncated (conditioned) to $(-\infty,x)$ for some $x\in\scrs$.  If $b(\cdot)$ is bounded away from both~$0$ and~$\infty$ on $\scrs$, then the constructions of the approximating primal and dual chains are identical to the case where 0 is reflecting, and details are omitted.  However, if $\lim_{x\rightarrow -\infty}b(x)=\infty$, then the
 approximating sequences of Markov chains need to be defined differently.  

To this end, on $\scrs^{n}:=\{S(1)-i_{n}\gD,\ldots, S(1)-\gD,S(1)\}$, with $i_{n}$ chosen so that $i_{n}\gD\rightarrow\infty$,
define a birth-and-death transition matrix $P^{n}$ via (here using the shorthand~$i$ for $S(1)-i\gD$)
$$P^{n}(i,i+1)=P^{n}(i,i-1):=\frac{b(i)h_{n}}{2\Delta^2}\text{ for }0<i<i_{n}$$
$$ P^{n}(i_{n},i_{n}-1):=\frac{b(i_{n})h_{n}}{\Delta^2},$$
$$P^{n}(0,1):=\frac{b(0)h_{n}}{\Delta^2},$$
with $P^{n}(i,i)$ chosen to make  the row sums of $P^{n}$ equal to 1, and
$$h_{n}:=\frac{\gD^2}{2\cdot \sup_{i\leq i_{n}}b(i)}$$
chosen again to ensure monotonicity.
For an initial probability distribution $\pi_0^{n}$ on $\scrs^{n}$, consider a birth-and-death Markov chain $X^{n}\sim (\pi_0^{n},P^{n})$.
Let the stationary distribution of $X^{n}$ be denoted $\pi^{n}$.  Let 
$$i_{n,x}:=\lfloor [S(1)-x]/\gD\rfloor,$$ 
and assume that $\pi_0^{n}:=\pi^{n,i_{n,x}}$ is $\pi^{n}$ truncated (conditioned) to 
$\{S(1)-i_n \gD,\ldots, S(1)-i_{n,x}\gD\}$.  Again note that $\pi_0^{n}\Rightarrow\pi^{(x)}$.

The following theorem is proven in a similar fashion to Theorem~\ref{thm:mccon}, and so the proof will be sketched with some detail omitted (see 
Appendix \ref{section:extpf} for notation).
\begt
\label{thm:entprimconv}
Assume $b(\cdot)$ is continuous and bounded away from~$0$ over $(-\infty,S(1)]$.  Let $\p Y^*(0)^{-1}=\pi^{(x)}$ for fixed $x\in\scrs\setminus\{-\infty\}$ and, as above, let $X^{n}\sim(\pi_0^{n},P^{n})$ with $\pi_0^{n}$ equal to $\pi^{n}$ truncated (conditioned) to  $\{S(1)-i_{n}\gD,\ldots, S(1)-i_{n, x}\gD\}$.  Define the continuous-time stochastic process $Y^{n}$ by setting 
$Y^{n}_t := X^{n}_{\lfloor t/h_{n}\rfloor}$ for $t \geq 0$. 
Then 
$Y^n \Rightarrow Y$.  
\ent
\begin{proof}
Fix $R$ such that $S(1)<R<\infty$.  With
$$
\tau_R^{n} := \inf\{t\geq 0:|Y^{n}_t| \geq R\},
$$
consider 
$$
Z^{n}(\cdot)
:=Y^{n}(\cdot\wedge \tau^{n}_R)
=X^{n}\bigg(\bigg\lfloor \frac{\cdot \wedge \tau_R^{n}}{h_{n}}\bigg\rfloor  \bigg).
$$  
With 
$$
\tau_R:=\inf\{t\geq 0:|Y_t|=R\},
$$
let $Z(\cdot):=Y(\cdot\wedge\tau_R)$.  
Denote the generator of~$Z$ by $A_Z$, with domain $\mathcal{D}(A_Z)$.
Writing $(T^{n}_R)$ for the transition semigroup associated with the Markov chain $X^{n}$ absorbed at any value $\leq -R$, let
$A^{n}_R:=h_{n}^{-1}(T^{n}_R-I)$.  Just as we showed~\eqref{eq:corem} from
Theorem \ref{thm:converge}, here we can show that 
\begin{equation}
\lim_{n \to \infty}\sup_{x\in[-R,S(1)]}\big| (A^{n}_R \rho_{n}f)(x)-(A_Z f)(x)\big|=0
\end{equation}
for all $f\in\mathcal{D}(A_Z)$.

By \cite[Corollary 4.8.9]{etku}, we have (see Appendix \ref{section:extpf}) that $Z^{n}\Rightarrow Z$.  The proof is finished by applying \cite[Theorem 11.1.1]{SV} to see that $Y^{n}\Rightarrow Y$ as desired.
\end{proof}

Let $Y^*$ be a SSD of $Y$, and let $Z^*$ be $Y^*$ put into natural scale.  
Form the dual Markov chain to $X^{n}$, and denote the dual by $\widehat{X}^{n}\sim(\delta_x,\widehat{P}^{n})$.  
Here $x\in\scrs\setminus\{-\infty\}$ 
is fixed and $Y_0\sim\pi^{(x)}$ (and hence $Y^*_0\sim\delta_x$). The following proposition gives the dual-convergence theorem analogous to Theorem \ref{thm:dualconv}.
\begt
\label{thm:entdualconv}
With the same assumptions as in Theorem \ref{thm:entprimconv}, further assume $b_Y(\cdot)\in C^2(-\infty,S(1)]$ and
\begin{equation}
\label{techcond}
\inf_{y\in\scrs} y^4 m_Y(M_Y^{-1}(-1/y))>0.
\end{equation}
For $t \geq 0$, define $\widehat Y^{n}(t) := \widehat X^{n}(\lfloor t/h_{n}\rfloor)$.
Then $\widehat Y^{n} \Rightarrow Y^*$.

\ent
\begin{proof}
The proof follows along the same path as the proof of Theorem \ref{thm:dualconv} and so details are omitted.  The only wrinkle here is the assumption~\eqref{techcond}, which is a technical condition needed to make the infinitesimal variance of the dual diffusion in natural scale bounded away from~$0$, which we exploited
in the proof of \refT{thm:dualconv}.
\end{proof}

\begr
Under some mild assumptions, the above theory can easily be extended to the case where both $0$ and $1$ are entrance boundaries for $X$.    For example, if~$X$ is in natural scale, it is sufficient that $b_X$ is bounded away from~$0$ and twice continuously differentiable on $\mathbb{R}$ and that the analogue of (\ref{techcond}) holds for $X$.  The analogues of  Theorem \ref{thm:mccon} and  Theorem \ref{thm:dualconv} can be easily recovered.  Details are omitted.
\enr

\section{Separation and hitting times}
\label{S:SST}
In 
the Markov chain setting, strong stationary duality gives that the separation distance in the primal chain is bounded by the tail probability of a suitable absorption time in the dual chain.  By studying and bounding the absorption time, which is sometimes more tractable than direct consideration of the separation distance, we can tightly bound the separation distance from stationarity in our primal chain.  See \cite{DFSST} for further detail.  Spelling this out more fully, if 
$X\sim(\pi_0,P)$ is an ergodic discrete-time Markov chain with state space $S$, stationary distribution $\pi$, and with SSD (as defined in \cite{DFSST}) $X^*\sim(\pi_0^*,P^*)$ absorbing in $m$, then for every $t$ we have the separation distance (which is a slight abuse of terminology, as separation is not a true distance)
\begin{equation}
\label{eq:dissep}
\sep(t):=\sup_{i\in S}\left(1-\frac{\pi_t(i)}{\pi(i)}\right)\leq \p_{\pi_0^*}(T_m^*>t),
\end{equation}
where $T_m^*$ is the hitting (i.e., absorption) time of state $m$ for the dual chain.
Under some monotonicity conditions, for example if the primal is a monotone likelihood ratio chain on a linearly ordered state space, the inequality in (\ref{eq:dissep}) can be made to be an equality for every~$t$ by a suitable formation of the dual chain.  

In 
our present diffusion setting, with $X$ a regular diffusion on $[0,1]$ with either reflecting or entrance behavior at the boundaries, we would like to recover a result similar to (\ref{eq:dissep}).  Let $\Pi$ be the invariant distribution for $X$, let  $X_0\sim \Pi_0$, and, given $t > 0$, let $\Pi_t$ be the corresponding distribution of $X_t$.  If $\Pi_t \ll \Pi$, 
define 
$$
a(t):=\essinf R_t = \sup \left\{r\,\big|\,\Pi(R_t < r) = 0 \right\}
$$ 
to be the essential infimum (with respect to~$\Pi$) of (any version of) the Radon--Nikodym derivative $R_t := d\Pi_t/d\Pi$.   We define the \textit{separation distance} (or \textit{separation}) of the diffusion from $\Pi$ at time $t$ as follows:
\begin{equation}
\label{sep}
\sep(\pi_t,\pi):=1-a(t).
\end{equation}
To simplify the notation, we shall write $\sep(t)$ for $\sep(\pi_t,\pi)$ unless the full notation is needed to avoid confusion.  

\begcl
Let $\sep(t)=\sep(\pi_t,\pi)$ be defined as above.  Then
\begin{enumerate}
\item[(a)] We have $0 \leq \sep(t) \leq 1$.
\item[(b)] For each~$t$ we have $\sep(t) = 0$ if and only if $\Pi_t = \Pi$.
\item[(c)] For any $\Pi_0$ we have $\Pi_t \ll \Pi$ for all $t > 0$.
\item[(d)] The separation $\sep(t)$ is non-increasing in~$t$.
\end{enumerate}
\encl
\begin{proof}
For~(a), we 
show equivalently that $0\leq a(t)\leq 1$.  To this end, 
let $R_t$ be (any version of) the Radon--Nikodym derivative $d\Pi_t/d\Pi$.   Since $R_t(y)\geq 0$ for all $y$, we have $a(t)\geq 0$.  But also
\begin{equation}
\label{eq:sepa}
1 = \int_0^1\!\Pi_t(dy) = \int_0^1\!R_t(y)\,\Pi(dy) \geq a(t)\int_0^1\!\Pi(dy) = a(t),
\end{equation}
finishing the proof.

For (b), note that if $\Pi_t=\Pi$, we can take $R_t \equiv1 $ as a version of the Radon--Nikodym derivative 
$d\Pi_t/d\Pi$, and thence $\sep(t) = 0$.  Conversely, if $\sep(t)=0$, then $a(t) = 1$ and \eqref{eq:sepa} is an equality; therefore $R_t = 1$ almost surely with respect to $\Pi$, and so $\Pi_t = \Pi$.

For (c), let $x \in (0,1)$.  When $\Pi_0=\delta_x$, regularity of $X$ guarantees the existence of a density for $\Pi_t$ with respect to $\Pi$, call it $f_x(\cdot)$.  For any $\Pi_0$, it follows that the $\Pi_0$-mixture of the densities $f_x(\cdot)$ is a density for $\Pi_t$ with respect to $\Pi$ [and so $\sep(t)$ is well defined].

For (d), for each $s>0$ let $R_{s} = d\Pi_s / d\Pi$.  Let $0<t<u$ and note for any 
$A \in \mathcal{B}$, the Borel $\sigma$-field of $[0,1]$, that
\begin{align*}
\int_A\!R_u(y)\,\Pi(dy) &= \Pi_u(A) = \int_0^1\!P_{u-t}(x,A)\,\Pi_t(dx) = \int_0^1\!R_t(x)\,P_{u-t}(x,A)\,\Pi(dx)\\
&\geq a(t) \int_0^1\!P_{u-t}(x,A)\,\Pi(dx) = a(t)\,\Pi(A)=\int_A\!a(t)\,\Pi(dy).
\end{align*}
Hence $R_u \geq a(t)$ almost surely with respect to~$\Pi$.  Hence $a(u)\geq a(t)$, and therefore $\sep(u) \leq\sep(t)$, as desired.
\end{proof}

As in the discrete setting, we are able to bound $\sep(t)$ in our primal diffusion~$X$ using the absorption time in state $1$ of our dual diffusion.  In the diffusion setting, by virtue of diffusions being stochastically monotone, the inequality in (\ref{eq:dissep}) is an equality without needing further assumptions.  Spelling this out:
\begl
\label{lem:sep}
Let~$X$ be 
a regular diffusion on $[0,1]$ begun in $\Pi_0$, let~$X$ have either reflecting or entrance behavior at the boundaries, and let $\Pi$ be the stationary measure for $X$.  Let $T^*_1$ be the hitting time of state 1 in the SSD diffusion  $X_t^*$ (as defined in Definition \ref{D:1.1}) begun in $\Pi_0^*$ satisfying \eqref{eq:init}.  Then
$$
\sep(t)=\p_{\Pi^*_0}(T^*_1>t)=1-\p_{\Pi_0^*}(X_t^*=1).
$$
\enl
\begin{proof}
Let $f \in \mathcal{F}[0,1]$.  By Remark \ref{rem:inter}, we have for all $t>0$ that $(\Pi_t,f)=(\Pi_t^*,\gL f)$.  Therefore, writing $R_t = d\Pi_t / d\Pi$ as usual,
we have 
\begin{align*}
\int_{[0,1]}\!\pi(x) R_t(x) f(x)\,dx
&= \int_{[0,1]}\!\Pi(dx)\,R_t(x) f(x)\\
&= \int_{[0,1]}\!\Pi_t(dx)\,f(x)\\
&= \int_{[0,1]}\!\Pi_t^*(dx)\int_{[0,x]}\!\pi^{(x)}(y)f(y)\,dy\\
&= \int_{[0,1]}\!\int_{[y,1]}\!\Pi^*_t(dx)\,\pi^{(x)}(y)f(y)\,dy.
\end{align*}
This holds for all $f\in\mathcal{F}[0,1]$, and so 
\begin{equation}
\label{sep3}
R_t(y)=\int_{[y,1]}\!\frac{\Pi^*_t(dx)}{\Pi(x)}
\end{equation}
for Lebesgue-a.e.\ (i.e.,\ for $\Pi$-a.e.)\ $y$.  Thus $\Pi(R_t < r) = 0$ if and only if the right side of~\eqref{sep3} is at least~$r$ for $\Pi$-a.e.~$y$, or, equivalently, $\Pi^*_t(\{1\})/\Pi(1)=\Pi^*_t(\{1\})\geq r.$  Therefore $a(t)=\Pi^*_t(\{1\})=\p_{\Pi_0^*}(X_t^*=1)$ and so $\sep(t)=1-\p_{\Pi_0^*}(X_t^*=1)$.
\end{proof}

\begin{remark}
We can also prove Lemma \ref{lem:sep} by passing to the limit the corresponding discrete-time results for the Markov chains in Section \ref{sec:2}.  First, suppose that $Y_0\sim \Pi^{(x)}$ for some $x> 0$ and hence $\Pi_0^*=\delta_x$ (see Remark \ref{rem:init}).
 Adopting the notation of Section \ref{sec:2}, the primal birth-and-death Markov chain $X^{n}\sim(\pi_0^{n}=\pi^{n,i_{n,x}},P^{n})$ has 
\begin{align*}
\sep^{n}(t)&=\sup_i\left(1-\frac{\sum_j \pi_0^{n}(j)P^{n}_t(j,i)}{\pi^{n}(i)}\right).
\end{align*}
Now
\begin{align*}
\frac{\sum_j \pi_0^{n}(j)P^{n}_t(j,i)}{\pi^{n}(i)}&=\frac{1}{H^{n}( i_{n,x})}\sum_{j\leq i_{n,x}}\frac{\pi^{n}(j)P^{n}_t(j,i)}{\pi^{n}(i)}\\
&=\frac{1}{H^{n}(i_{n,x})}\p_i(X^{n}_t\leq i_{n,x}).
\end{align*}
The monotonicity conditions outlined in Remark \ref{remark:conddual} and \cite[Remark 4.15]{DFSST} imply that this last expression
is minimized (for each $t=0,1,\ldots)$ when $i=n$, and that the minimum value is 
\begin{equation}
\frac{1}{H^{n}(i_{n,x})}\p_{n^{n}}(X^{n}_{t} \leq i_{n,x})=1-\sep^{n}(t)=\p_{i_{n,x}}(\widehat{T}_{n}\leq t),
\end{equation}
where $\widehat{X}^{n}$ is the strong stationary dual of $X^{n}$ as defined at 
\eqref{eq:death}--\eqref{eq:absorb}, with absorption time $\widehat{T}_{n}$ in its largest state $n$.  We now substitute $\lfloor t/h\rfloor$ for $t$, and recall that $h\equiv h_{n}$ is a function of $n$ and that $Y^{n}_t:=X^{n}_{\lfloor t/h \rfloor}$ (and analogously for $\widehat{Y}^{n}_t$), to find for real $t\geq0$ that
\begin{equation}
\label{eq:sepdt}
\frac{1}{H^{n}(i_{n,x})}\p_{S(1)}(Y^{n}_{t} \leq i_{n,x})=1-\sep^{n}(t)=\p_{i_{n,x}}\left(\widehat{Y}_{t}^{n}=S(1)\right). 
\end{equation}

With the necessary assumptions on $b_Y$ from \refT{thm:mccon} and \refT{thm:dualconv} or \refT{thm:entprimconv} and \refT{thm:entdualconv}, depending on the boundary assumptions of $X$, the left side of (\ref{eq:sepdt}) converges to $\frac{1}{\Pi(x)}\p_{S(1)}(Y_t\leq x)$ [where we note that the hypothesis of the theorems is met for the deterministic initial conditions $Y^{n}_0=Y_0=S(1)$].  \refT{thm:dualconv} implies that for any $\epsilon>0$ we have
\begin{equation}
\label{eps}
\lim_{n \to \infty} \p _{i_{n,x}}(\widehat Y_t^{n}> S(1)-\varepsilon)=\p_x(Y^*_t> S(1)-\varepsilon).\end{equation}
Let $\check{X}^{n}$ be the Siegmund dual of (the time-reversal of) $X^{n}$; by definition, 
$\check{X}^{n}$ is a Markov chain satisfying
$$\p_y(X^{n}_t\leq z)=\p_z(y\leq \check{X}^{n}_t)$$ for all $y,\ z\in\scrs^{n}$ and $t=0,1,2,\ldots$.  Equation~(5.3) in~\cite{DFSST} gives, with
$h = h_{n}$, $\gD=\gD_n$ and with $\lceil x\rceil_{\gD}$ (respectively, $\lfloor x\rfloor_{\gD}$) being
the smallest element $\geq x$ (resp.,\ the largest element $\leq x$) in the grid 
$\{S^*(S(0)),\,S^*(S(0)+\gD),\ldots, S^*(S(1)-\gD),\,S^*(S(1))\}$, that
\begin{align*}
\p _{i_{n,x}}(\widehat Y_t^{n}> S(1)-\varepsilon)&=\p _{i_{n,x}}(\widehat X_{\lfloor t/h\rfloor}^{n}> S(1)-\varepsilon)\\
&=\sum_{j>S(1)-\varepsilon} \frac{H^{n}(j)}{H^{n}(i_{n,x})}\p_{i_{n,x}}(\check X^{n}_{\lfloor t/h\rfloor}=j)\\
&\geq \frac{H^{n}(\lfloor S(1)-\varepsilon\rfloor_{\gD})}{H^{n}(i_{n,x})}\p_{i_{n,x}}(\check X^{n}_{\lfloor t/h\rfloor}> \lceil S(1)-\varepsilon\rceil_{\gD})\\
&\geq \frac{H^{n}(\lfloor S(1)-\varepsilon\rfloor_{\gD})}{H^{n}(i_{n,x})}\p_{i_{n,x}}(\check X^{n}_{\lfloor t/h\rfloor}\geq \lceil S(1)-2\varepsilon\rceil_{\gD})\\
&=\frac{H^{n}(\lfloor S(1)-\varepsilon\rfloor_{\gD})}{H^{n}(i_{n,x})}\p_{\lceil S(1)-2\varepsilon\rceil_{\gD}}( X^{n}_{\lfloor t/h\rfloor}\leq i_{n,x})\\
& \rightarrow \frac{\Pi(S(1)-\varepsilon)}{\Pi(x)}\p_{S(1)-2\varepsilon}(Y_t\leq x)\text{ as }n \to \infty,
\end{align*}
and this last expression converges to $\frac{1}{\Pi(x)}\p_{S(1)}(Y_t\leq x)$ as $\varepsilon \downarrow 0$.
We can also get an upper bound on $\p _{i_{n,x}}(\widehat Y_t^{n}> S(1)-\varepsilon)$ 
using
\begin{align*}
\sum_{j>S(1)-\varepsilon} \frac{H^{n}(j)}{H^{n}(i_{n,x})}\p_{i_{n,x}}(\check X^{n}_{\lfloor t/h\rfloor}=j)&\leq \frac{1}{H^{n}(i_{n,x})} \p_{i_{n,x}}(\check X^{n}_{\lfloor t/h\rfloor}\geq\lfloor S(1)-\varepsilon\rfloor_{\gD})\\
&=\frac{1}{H^{n}(i_{n,x})}\p_{\lfloor S(1)-\varepsilon\rfloor_{\gD}}( X^{n}_{\lfloor t/h\rfloor}\leq i_{n,x})\\
& \rightarrow \frac{1}{\Pi(x)}\p_{ S(1)-\varepsilon}(Y_t\leq x)\text{ as }n \to \infty,
\end{align*}
and this last expression converges to $\frac{1}{\Pi(x)}\p_{S(1)}(Y_t\leq x)$ as $\varepsilon \downarrow 0$.
Letting $\varepsilon \downarrow 0$ in (\ref{eps}) gives
$$\frac{1}{\Pi(x)}\p_{S(1)}(Y_t\leq x)=\p_x(T^*_{S(1)}\leq t).$$

Now $\Pi_t\ll \Pi$ for all $t>0$; let $R_t = 
d\Pi_t/d\Pi$, so that for any $A=[S(0),s)$ with $s\in\scrs$ we have
$$\p_{\Pi^{(x)}}(Y_t\in A)=\int_A \!R_t(y)\,\Pi(dy),$$
and also
\begin{align*}
\p_{\Pi^{(x)}}(Y_t\in A)&=\int_{[S(0),x]}\! \frac{\Pi(dy)}{\Pi(x)}P_t(y,A)\\
&=\int_{[S(0),x]} \frac{1}{\Pi(x)}\pi(y)\int_A\! p_t(y,z)\,dz\,dy.
\end{align*}
We will now appeal to the reversibility of $Y$.   A diffusion process $X$ with generator $A$ and state space $I$ is 
 \textit{reversible} with respect to the distribution $\mu$ if for all $f,\ g\in\mathcal{D}_A$ we have
\begin{equation}
\label{eq:rev}
\int\!f(y)\,(Ag)(y)\,\mu(dy)=\int\!(Af)(y)\,g(y)\,\mu(dy).
\end{equation}  
If $Y$ satisfies the assumptions of Lemma \ref{lem:sep}, noting that $f$, $g\in\mathcal{D}_A$ implies that the derivatives of each function vanish at the boundary of the state space, integration by parts yields that (\ref{eq:rev}) holds for $\mu=\Pi$, the stationary distribution of $Y$, and the primal diffusion is reversible with respect to $\Pi$.  Also note that (\ref{eq:rev}) is equivalent to the following (see \cite[Section II.5]{Lig}):  for all $f,g\in C(\scrs)$, and for all $t>0$ we have
\begin{equation}
\label{eq:rev2}
\int\!f(y)\,(T_tg)(y)\,\Pi(dy) = \int\!(T_tf)(y)\,g(y)\,\Pi(dy),
\end{equation}
where $(T_t)$ is the one parameter semigroup associated with $Y$.

Letting~$f$ and~$g$ be suitably continuous approximations of $\mathbbm{1}([S(0),x])$ and $\mathbbm{1}(A)$, and appealing to~\eqref{eq:rev2}, we have 
\begin{align*}
\int_{[S(0),x]} \frac{1}{\Pi(x)}\pi(y)\int_A\! p_t(y,z)\,dz\,dy
&=\int_{[S(0),x]} \frac{1}{\Pi(x)}\int_A p_t(z,y) \pi(z)\,dz\,dy\\
&=\int_A \frac{1}{\Pi(x)} \p_z(Y_t\leq x)\,\Pi(dz),
\end{align*}
and so $\frac{1}{\Pi(x)} \p_z(Y_t\leq x)$ is a version of $R_t(z)$.  By monotonicity of $Y$, we have that $\frac{1}{\Pi(x)} \p_z(Y_t\leq x)$ is minimized when $z=S(1)$, and hence for $Y$ we have  
$$a(t)=1-\sep(t)=\frac{1}{\Pi(x)} \p_{S(1)}(Y_t\leq x)=\p_x(T^*_{S(1)}\leq t),$$
establishing Lemma \ref{lem:sep} for $Y$.
\end{remark}

\begr
\label{rem:duality}
In the Markov chain setting of \cite{DFSST} and \cite{SSDCT}, the authors were able to justify their ``strong stationary duality'' nomenclature by tying their then-new notion of duality to the more classical notions of duality in the stochastic process literature.  Specifically, let $X\sim(\pi_0,P)$ be an
ergodic Markov chain with stationary distribution $\pi$.  If $X$ satisfies specific monotnicity conditions, namely, that the time reversal $\tilde{P}$ is monotone and $\pi_0(x) / \pi(x)$ decreases in $x$, then with $H$ be cumulative of~$\pi$, they show that the SSD $X^*$ of $X$ is the Doob H-transform of the Siegmund dual of the time-reversal of $X$.  

For a Markov process $X$ with transition operator $P_t(x,dy)$ and bounded 
harmonic function $H$ 
(satisfying $P^t H= H$ for all $t$), 
the \emph{Doob H-transform} of $X$ is the 
right-continuous Markov  process with transition operator
$$Q_t(x,dy):=\frac{H(y)}{H(x)}P_t(x,dy).$$
It has played a central role in Markov process duality theory, especially in the context of processes conditioned to die in a given set or point. See \cite[Chapter VII]{sharpe} for further detail.  The \emph{Siegmund dual} of a Markov process $X$ with ordered state space $\mathcal{X}$ is a Markov process $Y$ on $\mathcal{X}$ satisfying:
$$\p_x(X_t\leq y)=\p_y(x\leq Y_t)\text{ for all }x,y\in\mathcal{X}.$$
It has played a prominent role in the study of birth-and-death chains and one-dimensional diffusion theory and in the study of interacting particle systems (see \cite[Section II.3]{Lig} for extensive background).

To justify the nomenclature in the present diffusion setting, consider the diffusion $X$ as defined in Section \ref{S:Section 1}, and let $X^*$ be the strong stationary dual of $X$ specified in Definition~\ref{D:1.1}.  
Then, recalling from Remark~\ref{rem:sgmble} that for all $f\in \mathcal{F}[0,1]$ we have $\gL T_tf=T_t^*\gL f$, a simple calculation 
yields
$$
\int_{[0,x]}\!p_t(z,y)\,dy=\int_{[z,1]}\!\frac{\Pi(x)}{\Pi(y)}P^*_x(X^*_t\in dy),
$$ 
giving us immediately
that $X^*$ is the  Doob $H$-transform of the Siegmund dual  of (the time reversal of) $X$, where~$H$ here is the cumulative stationary distribution $\Pi$.  

A functional definition of duality generalizing Siegmund's definition was introduced in \cite{HS}.  For extensive background see again \cite[Section II.3]{Lig}.  Briefly, let~$X$ and~$Y$ be two Markov processes with ordered
state spaces $\mathcal{X}$ and $\mathcal{Y}$ and let $f$ be a bounded measurable function on $\mathcal{X}\times \mathcal{Y}$.  We define $Y$ to be the \emph{dual} of $X$ with respect to the function $f$ if 
$$\e_{x}f(X_t,y)=\e_{y}f(x,Y_t),\ \ \text{ for all }x\in\mathcal{X},\ y\in\mathcal{Y}.$$
As in \cite[Theorem 5.12]{DFSST}, a simple calculation yields that, in the diffusion setting, $X$ and its SSD $X^*$ are dual with respect to the function 
$$f(x,x^*):=\begin{cases}
1/\Pi(x^*),&\text{ if }x\leq x^*\\
0,&\text{ otherwise,}
\end{cases}$$
on $I\times I$, further justifying the duality name for $X^*$.  

With \cite[Definition 5.16]{DFSST}, the authors generalized the classical notion of functional duality from \cite{HS}.  Adapted to the present setting, let $X$ and $Y$ be two diffusions defined on a common probability space with state spaces $\mathcal{X}$ and $\mathcal{Y}$.  We say $Y$ is \emph{dual} to $X$ with respect to a function $f:\mathcal{X}\times\mathcal{Y}\rightarrow \mathbb{R}$ and distribution $\mu$ on $\mathcal{X}\times\mathcal{Y}$ if
$$\e_{\mu} f(X_t,Y_0)=\e_{\mu}f(X_0,Y_t).$$
In \cite[Theorem 5.19]{DFSST}, the authors were able to show that the strong stationary dual of an
ergodic Markov chain $X\sim(\pi_0,P)$ with stationary distribution $\pi$, and with the additional properties that the time reversal $\tilde{P}$ is monotone and $\pi_0(x) / \pi(x)$ decreases in $x$, is dual to the primal chain with respect to this new functional definition, for suitable choices of~$f$ and~$\mu$.  We are able to recover the analogue of their Theorem~5.19 here, as it is easy to see that $X^*$ (the strong stationary dual of $X$) and~$X$, appropriately coupled, are dual with respect to the function 
$f(x^*, x) = {\bf 1}(x \leq x^*) \pi(x) / \Pi(x^*)$ and $\mu$ equal to any mixture of the distributions $\delta_{x^*} \times \Pi^{(x^*)}$ with $x^* \in [0, 1]$.
\enr

\section{Hitting times and eigenvalues}
\label{S:HT}
In the continuous-time birth-and-death chain setting, a famous theorem 
due to
Karlin and MacGregor~\cite{kmac}
asserts
that the hitting time of state $n$ for a birth-and-death chain $X$ on $\{0,1,\ldots,n\}$ started in state $0$ is distributed as the sum of independent exponential random variables with parameters relating to the eigenvalues of the generator of $X$.  Fill~\cite{SSDCT} used strong stationary duality to exploit Karlin and MacGregor's result to prove that the separation from stationarity for an ergodic continuous-time birth-and-death chain $X$ at time $t$ is equal to $\p(Y>t)$ where $Y$ is a sum of independent exponential random variables with parameters depending on the eigenvalues of the generator of $X$.  In \cite{sepbd}, Diaconis and Saloff-Coste used Fill's result and tight concentration bounds on the tail probabilities of $Y$ to prove the existence of a separation cutoff for a sequence $(X_n)$ of birth-and-death chains under certain conditions on the eigenvalues of the generators of the chains $X_n$.  In this section, we outline and recover the analogous theory in the diffusion setting.

To this effect, consider again a diffusion $X$ on $[0,1]$  with generator $A$, and with reflecting or entrance boundary behavior at each boundary, satisfying the assumptions of Theorem~\ref{theorem:dual}.  Let $X^*$ be a strong stationary dual of $X$ according to Definition \ref{D:1.1}.  For fixed $\lambda$, let $v_{\lambda}(x)$ be the solution to the eigenvalue problem associated with $A$ (respectively, $A^*$):
\begin{align}
\label{eq:eval} Av+\gl v&=0\ \ \left(A^*v+\gl v=0\right)\text{ on }I^\circ
\end{align}
with boundary condition
\begin{align}
\label{eq:bound} B_0(v)=0
\end{align} where $B_0$ represents the following boundary condition (and $B^*_0$ the analogous dual boundary conditions at 0):
$$
B_0(v) :=
\begin{cases}
v(0), &\mbox{if $0$ is absorbing or exit;}\\
\frac{dv}{dS}^+(0), &\mbox{if $0$ is instantaneously reflecting or entrance}.
\end{cases}
$$
Let $T_{x,y}$ be the hitting time of $y$ for $X$ begun in $x$.  From \cite[Section 4.6]{ito}, we have that $v_{\lambda}(x)$ is unique up to multiplicative constant and that the moment generating function of $T_{x,y}$, call it $\psi_{x,y}$, can be expressed as 
\begin{equation}
\label{eq:mgf}
\psi_{x,y}(\lambda)=v_{\lambda}(x)/v_{\lambda}(y).
\end{equation}
A completely analogous set of results hold for $A^*$.

If we further add the relevant boundary condition at $1$, namely that $B_1(v)=0$ (where $B_1$ and $B_1^*$ are defined analogously to $B_0$), then we have from Sturm--Liouville theory (see for example \cite[Theorem 4.1]{kent}) that the eigenvalues of $-A^*$ (resp.,\ nonzero  
eigenvalues of $-A$) satisfying (\ref{eq:eval}) with the two boundary conditions are countable, real, positive, and simple and can be ordered such that
$$
0 < \gl_1 < \gl_2 < \cdots \uparrow\infty;
$$
further, they satisfy $\sum_{k=1}^{\infty}\lambda_k^{-1}<\infty$.  For extensive background on the relevant Sturm--Liouville theory, see for example~\cite{tit}.
The eigenfunctions and eigenvalues of $A$ and $A^*$ are connected by the following simple relationship:
\begp
\label{prop:evec}  Adopt the same assumptions as in Theorem \ref{theorem:dual}, 
and assume that $1$ is a reflecting boundary for $X$ and that $0$ is either a reflecting or entrance boundary for $X$.  Fix $\gl > 0$.
\smallskip 

$\mathrm{(a)}$  Suppose that $v = f$ is a solution of \eqref{eq:eval} for generator $A$ with boundary conditions $B_0(v)=B_1(v)=0$.  Then $v = \Lambda f$ is a solution of \eqref{eq:eval}  for generator $A^*$ with boundary conditions $B^*_0(v)=B^*_1(v)=0$ (and the same $\lambda$).
\smallskip

$\mathrm{(b)}$ 
Suppose that $v = g$ is a solution of \eqref{eq:eval} for generator $A^*$ with boundary conditions 
$B^*_0(v) = B^*_1(v) = 0$.  Then $f(\cdot)=g(\cdot)+\frac{\Pi(\cdot)}{\pi(\cdot)}g'(\cdot)$ is a solution 
of~\eqref{eq:eval}  for generator~$A$ with boundary conditions $B_0(v) = B_1(v) = 0$ (and the 
same~$\lambda$).
\enp

\begin{proof}
(a) If $f(\cdot)$ satisfies \eqref{eq:eval} for $A$ and the boundary conditions $B_0(f)=B_1(f)=0$, then 
$f\in\doma$ and 
$$
\frac{df}{dS}^+(0)=\left(\frac{f'}{s}\right)^+\!\!\!\!(0)
= 0 
=\frac{df}{dS}^-(1)=\left(\frac{f'}{s}\right)^-\!\!\!\!(1).
$$  
From~\eqref{eq:dom} we have $\Lambda f\in \domas$, and 
from~\eqref{eq:link} we have
\begin{equation}
\label{eq:dualevec}
A^* \Lambda f(\cdot)=\Lambda A f(\cdot)=\Lambda (-\lambda f)(\cdot)=-\lambda \Lambda f(\cdot)
\end{equation}
on $I^{\circ}$.  Therefore, $\Lambda f(\cdot)$ satisfies (\ref{eq:eval}) for $A^*$.  Also 
$B_0^*(\Lambda f)=0$ as 0 is an entrance boundary for the dual and $\Lambda f \in \domas$, and similarly $B_1^*(\Lambda f) = 0$.   

(b) Note that if $g(\cdot)$ satisfies~(\ref{eq:eval}) for $A^*$, then $g\in \domas$, and hence $g\in C[0,1]$, and $A^*g(1-)=-\lambda g(1-)=0$.  
Next, on $(0,1)$ note
\begin{align*}
f'&=g'+\frac{\pi^2-\Pi\pi'}{\pi^2}g'+\frac{\Pi}{\pi}g''=2g'-\frac{\Pi\pi'}{\pi^2}g'+\frac{\Pi}{\pi}g''\\
&=2g'+\frac{\Pi}{\pi}\left(g''+\frac{b'}{b}g'-\frac{2a}{b}g'\right)=2g'+\frac{\Pi}{\pi}\left(-\frac{2\lambda g}{b}-\frac{2\pi}{\Pi}g'  \right)\\
&=-2\lambda g \frac{\Pi}{\pi}\frac{1}{b}=-2\lambda g M s
\end{align*}
where the fourth equality follows from~\eqref{eq:eval}.  We have that $M^+(0) = 0 = g(1-)$ and 
$M(1-)<\infty$, and hence 
$$
\left(\frac{f'}{s}\right)^+(0) = 0 = \left(\frac{f'}{s}\right)^-(1),
$$ and therefore $B_0(f)=B_1(f)=0$.  
Next, on $(0,1)$ note
\begin{align*}
f''&=\frac{-2\lambda}{b}\bigg(\frac{\Pi}{\pi}g'+g\frac{\pi^2-\Pi\pi'}{\pi^2} -g\frac{\Pi b'}{\pi b}\bigg)\\
&=\frac{-2\lambda}{b}\bigg(\frac{\Pi}{\pi}g'+g+\frac{\Pi}{\pi}g \frac{s'}{s}\bigg)\in C(0,1),
\end{align*}
and hence $f\in C^2(0,1)$.
Combining the above, on $(0,1)$ we have
\begin{align*}
af'+\frac{1}{2}bf''&=\lambda g \frac{\Pi}{\pi}\frac{s'}{s}-\lambda\bigg(\frac{\Pi}{\pi}g'+g+\frac{\Pi}{\pi}g \frac{s'}{s}\bigg)\\
&=-\lambda\bigg(\frac{\Pi}{\pi}g'+g\bigg)=-\lambda f.
\end{align*}
To show that $f\in\doma$ and that $f$ satisfies (\ref{eq:eval}) for $A$ with the relevant boundary conditions it remains only to show that $f\in C[0,1]$.  We have (by \refT{theorem:dual}) that~$0$ is an entrance boundary for $X^*$, and hence for any fixed $\xi \in (0,1)$ we have that $N(0)<\infty$ and hence
$$
- \int_{(0, \xi]}\! S^*(\eta)\,M^*(d\eta) = \int_{(0, \xi]}\! \frac{1}{M(\eta)}s(\eta)M^2(\eta)\,d\eta
= \int_{(0, \xi]}\! \! s(\eta)M(\eta)\,d\eta < \infty,
$$
where the first equality holds by \eqref{eq:dualScale}--\eqref{eq:dualspeed}.  Since $g \in C[0, 1]$, it follows that
$$
\int_{(0, \xi]}\! |-2\lambda g(\eta)s(\eta)M(\eta)|\,d\eta =  \int_{(0, \xi]}\! |f'(\eta)|\,d\eta< \infty.
$$ 
Hence, by the dominated convergence theorem,
$$
\int_{(\omega, \xi]}\! f'(\eta)\,d\eta = f(\xi) - f(\omega)
$$
has a finite limit as $\omega \downarrow 0$.  We conclude that 
$f\in C[0,1)$.  We have by assumption that~$1$ is a reflecting boundary for $X$ and hence for any fixed $\xi \in (0,1)$ we have that $\Sigma(1)<\infty$ and hence $S[\xi,1)<\infty$ for all $\xi\in I^\circ$.  It follows that  
$$
\int_{[\xi, 1)}\!s(\eta)M(\eta)\,d\eta < \infty.
$$
By the same argument that showed that~$f$ is continuous at~$0$, we find that~$f$ is also continuous at~$1$. 
The proof is finished, as we have established that $f \in C[0,1]$.
\end{proof}

In the diffusions setting, we have an analogue (namely \cite[Theorem 5.1]{kent}) of Karlin and MacGregor's famous result on the eigenvalue expansion on birth-and-death hitting times.  Adapted to the present setting, we state the analogue as follows:
\begt
\label{theorem:eval} Let $X$ be a regular diffusion process on $[0,1]$ and assume $0$ is either instantaneously reflecting or entrance.  Then for $\gl<\gl_1$ we have
\begin{equation}
\label{eq:evalexp}
\lim_{x\rightarrow 0} \psi_{x,1}(\lambda)=\prod_{k=1}^{\infty}\bigg( 1-\frac{\lambda}{\lambda_k} \bigg)^{-1},
\end{equation}
which is the moment generating function of an infinite sum of independent exponential random variables with parameters $\lambda_k$.
\ent
Combining this with Proposition \ref{prop:evec} and Lemma \ref{lem:sep}, we arrive at
\begt
\label{T:sepeig}
Let~$X$ be a diffusion on $[0,1]$ with $X_0=0$, with generator $A$, and with either reflecting or entrance behavior at the bounday $0$ and reflecting behavior at the boundary $1$.   Let 
the eigenpairs 
$(\lambda_{i},v_{\lambda_{i}})$, $i = 1, 2, \dots$, of~$A$ with $\gl_i > 0$
satisfying~\eqref{eq:eval} and boundary conditions $B_0(v_{\lambda_i}) = 0 = B_1(v_{\lambda_i})$ be labeled so that $0 < \lambda_1 < \lambda_2 < \cdots$.  Let $X^*$ be a strong stationary dual of $X$ with generator 
$A^*$, and note that $X^*_0=0$ by Remark~\ref{rem:init}.  Let $W_1,W_2,\ldots$ be independent random variables with $W_i\sim \mbox{\rm Exp}(\lambda_i)$.  Then 
$$\sep(t)=\p_0(T_1^*>t)=\p(W>t)\ \text{where }\ W\stackrel{\mathcal{L}}{=}\sum_{i=1}^{\infty}W_i.$$ 
\ent
\noindent 
This mirrors the corresponding result for birth-and-death Markov chains given by \cite[Theorem 4.20]{DFSST} in discrete time and by \cite[Theorem 5]{SSDCT} in continuous time.  

In \cite{sepbd}, the authors used \cite[Theorem 4.20]{DFSST} to determine conditions for a separation cut-off to occur in a sequence of birth-and-death Markov chains.  We shall  presently derive analogous results for diffusions using Theorem \ref{T:sepeig}.
Consider now a sequence of diffusion generators $(A_n)_{n=1}^{\infty}$ defining a sequence of diffusions 
$(X^n)_{n=1}^{\infty}$  with $X^n_0\sim\nu^n,$ on compact intervals $[l_1,r_1]=I_1,\ [l_2,r_2]=I_2,\ \dots$ where all left boundary points, $l_n$ are assumed to be reflecting or entrance and all right boundary points 
$r_n$ are assumed to be reflecting.  Note that without loss of generality we can take $I_n=[0,r_n]$ for all $n\geq 1$.  We  write $\pi^n$ for the stationary distribution for $X^n$, and we write $\nu^n_t$ for 
the distribution of $X^n$ at time~$t$.  
This sequence of diffusions exhibits a \textit{separation cut-off at $(t_n)$} if the sequence $(t_n)$ is such that for any $\epsilon\in(0,1)$ we have
\begin{align*}
\mathrm{(i)}\ \lim_{n\rightarrow \infty} \sep(\nu^n_{(1+\varepsilon)t_n}, \pi^n)&=0, \text{ and } \\
\mathrm{(ii)}\ \lim_{n\rightarrow \infty} \sep(\nu^n_{(1-\varepsilon)t_n}, \pi^n)&=1.
\end{align*}

To apply Theorem \ref{T:sepeig} here, let the nonzero eigenvalues of $A_n$ be labeled $0<\lambda_{n,1}<\lambda_{n,2}<\cdots,$ and let $\nu^n=\delta_0$ for all $n\geq 1$.  We further assume that each $A_n$ satisfies the assumptions of Theorem \ref{thm:entrefl}, and let $(A_n^*)_{n=1}^{\infty}$ be the sequence of  generators of the strong stationary duals of $(X^n)_{n=1}^{\infty}$ as defined by Definition \ref{D:1.1}
For each $n\geq 1$, let $W_{n,j} \sim \mbox{\rm Exp}(\lambda_{n,j})$ be independent random variables, and let $W_n\stackrel{\mathcal{L}}{=}\sum_{j=1}^{\infty}W_{n,j}$.  From Theorem \ref{T:sepeig}, we have $\sep^n(t)=\p(W_n>t)$.  We can therefore get sharp bounds on separation by deriving sharp bounds for the tail probabilities of $W_n$.  To this end, note that
we have
$$
\e\,W_n = \sum_{j=1}^{\infty} \lambda_{n,j}^{-1}<\infty, \quad \va\,W_n = \sum_{j=1}^{\infty} \lambda_{n,j}^{-2}<\infty.
$$
An application of the one-sided Chebyshev's inequality gives the analogue to the separation cut-off result \cite[Theorem 5.1]{sepbd}: 
\begt
\label{th:sepcut}
Let $(A_n)_{n=1}^{\infty}$ be a sequence of diffusion generators defining diffusions $(X^n)_{n=1}^{\infty}$, with $X^n_0 \sim \nu^n=\delta_0$, on compact intervals $[0,r_1]=I_1,\ [0,r_2]=I_2,\ \dots$, where $0$ is assumed to be reflecting or entrance for all $n$, and all right boundary points $r_n$ are assumed to be reflecting. With the eigenvalues $\lambda_{n,i}$ defined as above,  this sequence of diffusions exhibits a separation cut-off 
if and only if
$$
\lim_{n\rightarrow\infty}\lambda_{n,1}\,\e W_n = \infty,
$$ 
in which case there is a separation cut-off at $(t_n)$ with $t_n:=\e W_n$.  
Further, for any $c > 0$ the following separation bounds hold for \emph{any} sequence $(t_n)$, 
where we restrict to $c\leq 1$ in the second bound:
$$\sep(\nu^n_{(1+c)t_n},\pi^n)\leq\frac{1}{1+c^2 \lambda_{n,1}t_n},\ \ \sep(\nu^n_{(1-c)t_n},\pi^n)\geq1-\frac{1}{1+c^2 \lambda_{n,1}t_n}.$$
\ent
\noindent The proof is completely analogous to the proof of Theorem 5.1 in \cite{sepbd}, and so is omitted.  

\begx
\emph{Let $0<1 = r_1 \leq r_2\leq r_3\leq\cdots$ be an arbitrary increasing sequence of positive real numbers, and let $A_n$ be the generator of reflecting Brownian motion on $I_n=[0,r_n]$.  Then (see \cite[Section 6]{kent}), we know that $$\lambda_{n,k}=\frac{j^2_{k}}{2 r_n^2}$$ where  
$j_k = k \pi$ is well known.  Note 
$$
\lambda_{n,1} \e W_n = \sum_{k=1}^{\infty}\frac{j^2_1}{j^2_k} = \frac{\pi^2}{6}
$$ 
is constant in $n$, and therefore there is no separation 
cut-off.
}
\enx

\begx
\emph{See Example \ref{ex:bes}.
Let $(\eta_n)$ be a sequence of positive real numbers diverging monotonically to infinity.  Let $A_n$ be the generator for a Bes(2$\eta_n$+2) process on $[0,1]$ with~$1$ a reflecting boundary.  Again from \cite[Section 6]{kent}, we have that 
$$\lambda_{n,k}=\frac{j^2_{n,k}}{2 }$$ where $(j_{n,k})_{k=1}^{\infty}$ are the positive zeros of the Bessel function $J_{\eta_n+1}$.  Then
$$\lambda_{n,1} \e W_n = \sum_{k=1}^{\infty}\frac{j^2_{n,1}}{j^2_{n,k}}.$$
It 
is well known (see for instance \cite[equations~(1) and~(40)]{sned}) that
$$\sum_{k=1}^{\infty}\frac{1}{j^2_{n,k}}=\frac{1}{4(\eta_n+2)}$$
and (see \cite[pg. 371]{absteg}) that 
$$j^2_{n,1}=\left[\eta_n+1+O\left(\eta_n^{1/3}\right)\right]^2;$$ so there is a separation cut-off for this sequence of diffusions at $(t_n)$, with $t_n =2\sum_{k=1}^{\infty}j_{n,k}^{-2}=(2\eta_n+4)^{-1}$.
}

\emph{This is, perhaps, not a surprising result in light of the interpretation of the Bes($m$) process as the radial part of $m$-dimensional Brownian motion for integer~$m$.  As the strong stationary dual of a Bes($\alpha$) process is a Bes($\alpha+2$) process, for integer sequences $\eta_n=m_n$, a separation cut-off is equivalent to a sharp concentration in the hitting time of $1$ of the dual Bes($2m_n+4$) sequence, i.e., a sharp concentration in the hitting time of the unit sphere for $(2m_n+4)$-dimensional Brownian motion started in $\vec{0}$. 
For large $m_n$, at time $t$ the ratio of the square of the radial part of $(2m_n+4)$-dimensional Brownian motion to $t$ has a distribution which doesn't depend on $t$ and (by the central limit theorem) is approximately normal with mean $2m_n+4$ and variance $2(2m_n+4)$.  We therefore expect to have a sharp concentration of the hitting time of the unit sphere at $t=(2m_n+4)^{-1}$, and indeed we found that the 
cut-off occurs there.
}
\enx

\bibliographystyle{plain}
\bibliography{05references}

\begin{thebibliography}{10}

\bibitem{absteg}
M.~Abramowitz and I.~A. Stegun.
\newblock {\em Handbook of Mathematical Functions with Formulas, Graphs, and
  Mathematical Tables}.
\newblock Dover, New York, ninth {D}over printing, tenth {GPO} printing
  edition, 1964.

\bibitem{bhat}
R.~N. Bhattacharya and E.~C. Waymire.
\newblock {\em Stochastic processes with applications}.
\newblock Wiley Series in Probability and Mathematical Statistics: Applied
  Probability and Statistics. John Wiley \& Sons Inc., New York, 1990.
\newblock A Wiley-Interscience Publication.

\bibitem{billing}
P.~Billingsley.
\newblock {\em Convergence of probability measures}.
\newblock Wiley Series in Probability and Statistics: Probability and
  Statistics. John Wiley \& Sons Inc., New York, second edition, 1999.
\newblock A Wiley-Interscience Publication.

\bibitem{DFSST}
P.~Diaconis and J.~A. Fill.
\newblock Strong stationary times via a new form of duality.
\newblock {\em Ann. Probab.}, 18(4):1483--1522, 1990.

\bibitem{PDLmBD}
P.~Diaconis and L.~Miclo.
\newblock On times to quasi-stationarity for birth and death processes.
\newblock {\em J. Theoret. Probab.}, 22(3):558--586, 2009.

\bibitem{sepbd}
P.~Diaconis and L.~Saloff-Coste.
\newblock Separation cut-offs for birth and death chains.
\newblock {\em Ann. Appl. Probab.}, 16(4):2098--2122, 2006.

\bibitem{dynk}
E.~B. Dynkin.
\newblock {\em Markov processes. {V}ols. {I}, {II}}, volume 122 of {\em
  Translated with the authorization and assistance of the author by J. Fabius,
  V. Greenberg, A. Maitra, G. Majone. Die Grundlehren der Mathematischen
  Wissenschaften, B\"ande 121}.
\newblock Academic Press Inc., Publishers, New York, 1965.

\bibitem{etku}
S.~N. Ethier and T.~G. Kurtz.
\newblock {\em Markov processes}.
\newblock Wiley Series in Probability and Mathematical Statistics: Probability
  and Mathematical Statistics. John Wiley \& Sons Inc., New York, 1986.
\newblock Characterization and convergence.

\bibitem{SSDCT}
J.~A. Fill.
\newblock Strong stationary duality for continuous-time {M}arkov chains. {P}art
  {I}: {T}heory.
\newblock {\em Journal of Theoretical Probability}, 5:45--70, 1992.
\newblock 10.1007/BF01046778.

\bibitem{fps}
J.~A. Fill.
\newblock An interruptible algorithm for perfect sampling via {M}arkov chains.
\newblock {\em Ann. Appl. Probab.}, 8(1):131--162, 1998.

\bibitem{JAFBD}
J.~A. Fill.
\newblock The passage time distribution for a birth-and-death chain: strong
  stationary duality gives a first stochastic proof.
\newblock {\em J. Theoret. Probab.}, 22(3):543--557, 2009.

\bibitem{fk}
J.~A. Fill and J.~Kahn.
\newblock Comparison inequalities and fastest-mixing markov chains.
\newblock {\em Annals of Applied Probability}, 23(5):1778--1816, 2013.

\bibitem{LF}
J.~A. Fill and V.~Lyzinski.
\newblock Hitting times and interlacing eigenvalues: A stochastic approach
  using intertwinings.
\newblock {\em Journal of Theoretical Probability}, pages 1--28, 2012.

\bibitem{MR2001m:60164}
J.~A. Fill, M.~Machida, D.~J. Murdoch, and J.~S. Rosenthal.
\newblock Extension of {F}ill's perfect rejection sampling algorithm to general
  chains.
\newblock {\em Random Structures Algorithms}, 17(3-4):290--316, 2000.
\newblock Special issue:\ Proceedings of the Ninth International Conference
  ``Random Structures and Algorithms'' (Poznan, 1999).

\bibitem{HS}
R.~Holley and D.~Stroock.
\newblock Dual processes and their application to infinite interacting systems.
\newblock {\em Adv. in Math.}, 32(2):149--174, 1979.

\bibitem{ito}
K.~It{\^o} and H.~P. McKean, Jr.
\newblock {\em Diffusion processes and their sample paths}.
\newblock Springer-Verlag, Berlin, 1974.
\newblock Second printing, corrected, Die Grundlehren der mathematischen
  Wissenschaften, Band 125.

\bibitem{kmac}
S.~Karlin and J.~McGregor.
\newblock Coincidence properties of birth and death properties.
\newblock {\em Pacific J. Math.}, 9:1109--1140, 1959.

\bibitem{kt2}
S.~Karlin and H.~M. Taylor.
\newblock {\em A second course in stochastic processes}.
\newblock Academic Press Inc. [Harcourt Brace Jovanovich Publishers], New York,
  1981.

\bibitem{kent}
J.~T. Kent.
\newblock Eigenvalue expansions for diffusion hitting times.
\newblock {\em Z. Wahrsch. Verw. Gebiete}, 52(3):309--319, 1980.

\bibitem{knight}
F.~B. Knight.
\newblock {\em Essentials of {B}rownian motion and diffusion}, volume~18 of
  {\em Mathematical Surveys}.
\newblock American Mathematical Society, Providence, R.I., 1981.

\bibitem{lamp}
J.~Lamperti.
\newblock {\em Stochastic processes}.
\newblock Springer-Verlag, New York, 1977.
\newblock A survey of the mathematical theory, Applied Mathematical Sciences,
  Vol. 23.

\bibitem{Lig}
T.~M. Liggett.
\newblock {\em Interacting Particle Systems}.
\newblock Springer Berlin Heidelberg, 1985.

\bibitem{VL}
V.~Lyzinski.
\newblock {\em Intertwinings, interlacing eigenvalues, and strong stationary
  duality for diffusions}.
\newblock PhD thesis, Johns Hopkins University, January 2013.

\bibitem{mic}
L.~Miclo.
\newblock Strong stationary times for one-dimensional diffusions.
\newblock {\em arXiv preprint arXiv:1311.6442}, 2013.

\bibitem{mish}
S.~Pal and M.~Shkolnikov.
\newblock Intertwining diffusions and wave equations.
\newblock {\em arXiv preprint arXiv:1306.0857}, 2013.

\bibitem{ROGPIT}
L.~C.~G. Rogers and J.~W. Pitman.
\newblock Markov functions.
\newblock {\em Ann. Probab.}, 9(4):573--582, 1981.

\bibitem{rogwil}
L.~C.~G. Rogers and D.~Williams.
\newblock {\em Diffusions, {M}arkov processes, and martingales. {V}ol. 2}.
\newblock Cambridge Mathematical Library. Cambridge University Press,
  Cambridge, 2000.
\newblock It{\^o} calculus, Reprint of the second (1994) edition.

\bibitem{sharpe}
M.~Sharpe.
\newblock {\em General theory of {M}arkov processes}, volume 133 of {\em Pure
  and Applied Mathematics}.
\newblock Academic Press Inc., Boston, MA, 1988.

\bibitem{sned}
I.~N. Sneddon.
\newblock On some infinite series involving the zeros of {B}essel functions of
  the first kind.
\newblock {\em Proc. Glasgow Math. Assoc.}, 4:144--156 (1960), 1960.

\bibitem{SV}
D.~W. Stroock and S.~R.~S. Varadhan.
\newblock {\em Multidimensional diffusion processes}, volume 233 of {\em
  Grundlehren der Mathematischen Wissenschaften [Fundamental Principles of
  Mathematical Sciences]}.
\newblock Springer-Verlag, Berlin, 1979.
\newblock Reprinted in 2006.

\bibitem{tit}
E.~C. Titchmarsh.
\newblock {\em Eigenfunction expansions associated with second-order
  differential equations. {P}art {I}}.
\newblock Second Edition. Clarendon Press, Oxford, 1962.

\end{thebibliography}
\appendix
\section{Details of Theorem \ref{thm:converge}}
\label{section:extpf}
In proving Theorem \ref{thm:mccon}, we made use of \refT{thm:converge} adapted from \cite[Theorems 4.8.2 and 1.6.5 and Corollary 4.8.9]{etku}].  We restate the theorem here:\\
\noindent{\bf Theorem 4.2}\emph{ Let~$A$ be the generator (as in Section \ref{S:Section 1}) of a regular diffusion process~$Y$ with state space~$\mathcal{Y}$.  Assume $h_n > 0$ converges to~$0$ as  $n \to \infty$.  Let $X^{n}\sim(\pi_0^{n},P^{n})$ be a Markov chain on metric state space $\mathcal{Y}^{n} \subset \mathcal{Y}$ and define $Y^{n}_t := X^{n}_{\lfloor t/h_{n} \rfloor}$.  Further assume $Y^{n}_0\Rightarrow Y_0$.  Letting $\mathcal{B}(\scry^n)$ be the space of real-valued bounded measurable functions on $\scry^n$, define $T^{n}:\mathcal{B}(\scry^{n})\rightarrow \mathcal{B}(\scry^n)$ via
$$
T^{n} f(x)=\e_x f(X^{n}_1).
$$ 
Let $\rho_{n}:C(\scry)\rightarrow \scrb(\scry^{n})$ be defined via $\rho_{n}f(\cdot)=f |_{\scry^{n}}(\cdot)$.  If $\doma$ is an algebra that strongly separates points, and
\begin{equation}
\label{eq:core}
\lim_{n\rightarrow \infty} \sup_{y\in\scry^{n} }\left| (A^{n}\rho_{n} f)(y) - (A f)(y) \right|=0
\end{equation}
for 
all $f\in \doma$, then $Y^n \Rightarrow Y$.}
\smallskip

The purpose of this appendix is to carefully spell out the proof of the above theorem, as the notation in \cite{etku} differs considerably from the notation we have adopted.  The following chart gives the notational equivalences between the present  work 
and~\cite{etku}; in connection with $\mu_n(x, \cdot)$, see Corollary~4.8.5 in~\cite{etku}.
\[ \begin{array}{l|r}
\mbox{Notation in present work:} & \mbox{Notation in \cite{etku}:}\\
\hline
\mathcal{Y}\mbox{ with the Euclidean metric} & (E, r) \\
\scry^{n},\ A^{n},\ T^{n} & E_n,\ A_n,\ T_n\\
P^{n}(x,\cdot) & \mu_n(x,\cdot)\\
\{(f,Af)\,|\,f\in\doma\} & A \\
\doma & \doma \\
C(\scrs) & L \\
1 / h_n & \alpha_n \\
\text{id} & \eta_n\\
\rho_{n} & \pi_n\\
\end{array}\]
Here $\rho_{n}:C(\scry)\rightarrow \mathcal{B}(\scry^{n})$ is defined via 
$\rho_{n}f(\cdot)=f |_{\scry^{n}}(\cdot)$, and $\text{id}:\scry^{n}\rightarrow \scry$ is the inclusion function embedding $\scry^{n}$ into $\scry$.

\begin{proof}[Proof of \refT{thm:converge}]
Clearly $C(\scry)$ is convergence determining, and by considering suitably smooth uniform approximations in $\doma$ to 
the indicator function of $\{x\}$
for each $x\in\scry$, it follows that $\doma\subset C(\scry)$ is an algebra that strongly separates points.  In the notation of \cite[Corollary 4.8.9]{etku}, we have $G_n=E_n=\scry^n$ , and so to 
prove $Y^n\Rightarrow Y$, it suffices to prove that for each $T > 0$ and $f \in C(\scry)$ we have
\begin{equation*}
\label{eq:semigrpconv}
\tag{B.2}
\lim_{n \to \infty}\sup_{y\in \scry^{n}}\left| (T^{n})^{\lf t / h \rf} \rho_{n} f(y)-\rho_{n} T_t f(y)\right|, \qquad 0\leq t\leq T.
\end{equation*}

From \cite[Theorem 1.6.5]{etku}, to prove (\ref{eq:semigrpconv}) it suffices to establish that for all 
$f\in\doma$ we have that 
$\rho_{n} f\in  \mathcal{B}(\scry^{n}) (=L_n$ in the notation of \cite[Theorem 1.6.5]{etku}) satisfies
\begin{equation*}
\label{eq:f}
\tag{B.3}
\lim_{n \to \infty} \sup_{y\in\scry^{n} }\left|\rho_{n} f(y) - f(y) \right| = 0
\end{equation*}
and 
\begin{equation*}
\label{eq:A}
\tag{B.4}
\lim_{n \to \infty} \sup_{y\in\scry^{n} }\left| (A^{n}\rho_{n} f)(y) - (A f)(y) \right|=0.
\end{equation*}
But (\ref{eq:f}) is clearly true, and~(\ref{eq:A}) is assumed (for all $f\in\doma$) in the statement of\ \refT{thm:converge}.
\end{proof}

\end{document}